\documentclass[10pt,a4paper]{article}
\usepackage{graphicx} % figures
\usepackage{amstext} % math
\usepackage[top=30mm,bottom=30mm,left=25mm,right=25mm]{geometry} % Two-column format
\usepackage{hyperref} % links and citations
\usepackage{amsmath, amsfonts, amssymb, amsthm, mathtools}
\usepackage{tikz}
\usepackage{listofitems}

\usetikzlibrary{decorations.pathmorphing,patterns}

\tikzstyle{mynode}=[thick,draw=blue,fill=blue!20,circle,minimum size=12,inner sep=0pt]
\tikzstyle{mynode2}=[thick,draw=red,fill=red!20,circle,minimum size=12,inner sep=0pt]

\newtheorem{remark}{Remark}

\begin{document}
% \title{\vspace{-18mm}
% \begin{minipage}{\linewidth}
% \hspace{5mm}\raisebox{-50pt}{\includegraphics[width=.23\textwidth]{logo-3AF}}\hspace{4mm}
% \textcolor{myblu}{\textbf{\textit{\normalsize\begin{tabular}{l}
% 58$^\text{th}$ 3AF International Conference\\ on Applied Aerodynamics\\ 27 --- 29 March 2024, Orléans -- France
% \end{tabular}}}}
% \hspace{12mm}\textbf{\normalsize AERO2024-07-DUVIGNEAU}
% \end{minipage}\\\vspace{10mm}
% % Paper title
% \textbf{\Large INVESTIGATIONS ON PHYSICS-INFORMED\\ NEURAL NETWORKS FOR AERODYNAMICS}}
% % Author(s): 10 pt bold (Author)
% \author{\textbf{\normalsize Guillaume Coulaud$^\text{(1)}$, Maxime Le$^\text{(1)}$ and Régis Duvigneau$^\text{(1)}$}
% % Affiliation(s): 10 pt italic (Affiliation)
% \\{\normalsize\itshape
% $^\text{(1)}$Université Côte d'Azur, Inria, CNRS, LJAD, 2004 route des lucioles 06902 Sophia-Antipolis, France}
% \\{\normalsize\itshape
% guillaume.coulaud@inria.fr, maxime.le@inria.fr, regis.duvigneau@inria.fr}
% }\date{}

\title{Investigations on Physics-Informed Neural Networks\\ for Aerodynamics}
\author{Guillaume Coulaud, Maxime Le \& Régis Duvigneau\\
\small{Université Côte d'Azur, Inria, CNRS, LJAD}\\
\small{2004 route des lucioles 06902 Sophia-Antipolis, France}}
\date{March 2024}

\maketitle
\begin{abstract}
Physics-Informed Neural Networks (PINNs) have recently emerged as a novel approach to simulate complex physical systems on the basis of both data observations and physical models. In this work, we investigate the use of PINNs for various applications in aerodynamics and we explain how to leverage their specific formulation to perform some tasks effectively. In particular, we demonstrate the ability of PINNs to construct parametric surrogate models, to achieve multiphysic couplings and to infer turbulence characteristics via data assimilation. The robustness and accuracy of the PINNs approach are analysed, then current issues and challenges are discussed.

\end{abstract}

\section{Introduction}

For decades, simulation methods in engineering rely on Partial Differential Equation (PDE) solvers, such as Finite-Volume or Finite-Element methods, which have proved their robustness and accuracy for most application domains. Recently, data-based approaches have emerged as possible concurrents, especially in problems for which PDE models are not well established, e.g. turbulence phenomena. However, in several engineering domains, data are not so easy to collect or generate, which reduces the range of applicability of this methodology. Physics-Informed Neural Networks (PINNs)~\cite{raissiPhysicsinformedNeuralNetworks2019} offer a compromise by dealing with both PDEs and data observations to construct a mixed model. Numerous articles have since been published in the literature to demonstrate their potentiality for different applications~\cite{cuomoScientificMachineLearning2022}.

In this work, we investigate the use of PINNs for non-classical problems, for which their specfic formulation can be leveraged to define an efficient resolution procedure. In particular, we investigate the construction of parametric models that allow to get an instantaneous response for a range of physical parameters, after the training has been achieved. Then, the use of PINNs to facilitate multidisciplinary couplings is studied. Finally, their ability to handle both models and data is exploited to devise an efficient data assimilation approach. In all cases, flow problems are used as illustrations.

\section{Methodology}

\subsection{PINNs principle}

The objective of PINNs is to simulate a physical system using a Multi-Layer Perceptron (MLP) neural network, which is trained according to known physical rules as well as a set of observation data~\cite{raissiPhysicsinformedNeuralNetworks2019}. For the sake of simplicity, we consider a generic problem to present the methodology, for which the physical rules are expressed as a first-order PDE, associated to Dirichlet and Neumann boundary conditions, and initial condition. The extension to more complex governing equations will be straightforward. Thus, the problem writes:
\begin{align}
\begin{cases}
    \frac{\partial{u}}{\partial{t}} + \mathcal{N}(u,\frac{\partial{u}}{\partial{x}}) = 0 & (x,t) \in \Omega \times [0,T] \\
    u(x^D,t) = g^D(x^D,t) & (x^D,t) \in \partial \Omega^D \times [0,T]\\
    \frac{\partial{u}}{\partial{n}}(x^N,t) = g^N(x^N,t) & (x^N,t) \in \partial \Omega^N \times [0,T]\\
    u(x,0) = g^I(x) & x \in \Omega \\
    u(x_i,t_i) = u^\star_i            & i=1,\ldots,N_{data} \\
\end{cases}
\label{eq:pb}
\end{align}

where $\Omega$ is an open domain, $(x,t)$ the space-time coordinates, $u$ the solution field and $\mathcal{N}$ denotes a first-order differential operator. $g^D$, $g^N$ and $g^I$ define respectively the Dirichlet, Neumann and initial conditions. Additionally, a set of observation data $(u^\star_1,\ldots,u^\star_{N_{data}})$ are provided.

\begin{figure}[!ht]
\centering
\begin{tikzpicture}[x=1cm,y=0.5cm]

\node at (4,-4.5)  {descent step};
\draw[->,thick] (6,-5) -- (3,-5) -- (3,-3);

\draw[-,thick] (6,-5) -- (11,-5);
\node[mynode2] at (6,-5) {$\frac{\partial}{\partial{\theta}}$};

\draw[-,thick] (15.5,2) -- (16,0);
\draw[-,thick] (14,0) -- (16,0);
\draw[-,thick] (14,-1) -- (16,0);
\draw[-,thick] (14,-2) -- (16,0);
\draw[-,thick] (14,-3) -- (16,0);
\draw[->,thick] (16,0) -- (16,-4);
\node at (14,-5)  {$\mathcal{L} = \mathcal{L}_{PDE} + \mathcal{L}_{I} + \mathcal{L}_{D} + \mathcal{L}_{N} + \mathcal{L}_{data}$} ;

\draw[-,thick] (6,-2) -- (7,-3);
\draw[-,thick] (5,0) -- (7,-2);
\draw[-,thick] (5,0) -- (7,-1);
\draw[-,thick] (5,0) -- (7,0);
\node[right] at (7,0)  {$\mathcal{L}_{data} = \frac{1}{N_{data}} \sum_{i=1}^{N_{data}} \left [ u(x_i,t_i) - u^\star_i \right ]^2$};
\node[right] at (7,-1)  {$\mathcal{L}_{I} = \frac{1}{N_{I}} \sum_{i=1}^{N_{I}} \left [ u(x_i,0) - g^I(x_i,0) \right ]^2$};
\node[right] at (7,-2)  {$\mathcal{L}_{D} = \frac{1}{N_{D}} \sum_{i=1}^{N_{D}} \left [ u(x_i^D,t_i) - g^D(x_i^D,t_i) \right ]^2$};
\node[right] at (7,-3)  {$\mathcal{L}_{N} = \frac{1}{N_{N}} \sum_{i=1}^{N_{N}} \left [ \frac{\partial{u}}{\partial{n}}(x_i^N,t_i) - g^N(x_i^N,t_i) \right ]^2$};

\draw[-,thick] (6,3) -- (7,2);
\draw[-,thick] (6,1) -- (7,2);
\node[right] at (7,2)  {$\mathcal{L}_{PDE} = \frac{1}{N_{PDE}} \sum_{i=1}^{N_{PDE}} \left [ \frac{\partial{u}}{\partial{t}}(x_i,t_i) + \mathcal{N}(u,\frac{\partial{u}}{\partial{x}})(x_i,t_i) \right]^2$};

\draw[-,thick] (5,0) -- (6,3);
\draw[-,thick] (5,0) -- (6,1);
\draw[-,thick] (5,0) -- (6,-2);
\node[mynode2] at (6,3)  {$\frac{\partial}{\partial{x}}$};
\node[mynode2] at (6,1)  {$\frac{\partial}{\partial{t}}$};
\node[mynode2] at (6,-2)  {$\frac{\partial}{\partial{n}}$};

\node at (2.5,-2.5)  {$\theta$};
  \readlist\Nnod{2,5,5,5,1} % number of nodes per layer
  \foreachitem \N \in \Nnod{ % loop over layers
    \foreach \i [evaluate={\x=\Ncnt; \y=\N/2-\i+0.5; \prev=int(\Ncnt-1);}] in {1,...,\N}{ % loop over nodes
      \node[mynode] (N\Ncnt-\i) at (\x,\y) {};
      \ifnum\Ncnt>1 % connect to previous layer
        \foreach \j in {1,...,\Nnod[\prev]}{ % loop over nodes in previous layer
          \draw[-,thick] (N\prev-\j) -- (N\Ncnt-\i); % connect arrows directly
        }
      \fi % else: nothing to connect first layer
    }
  }
\node at (1,0.5)  {$x$};
\node at (1,-0.5)  {$t$};
\node at (5,0)  {$u$};

\end{tikzpicture}
\caption{Overview of PINNs principle}
\label{fig:pinns}
\end{figure}
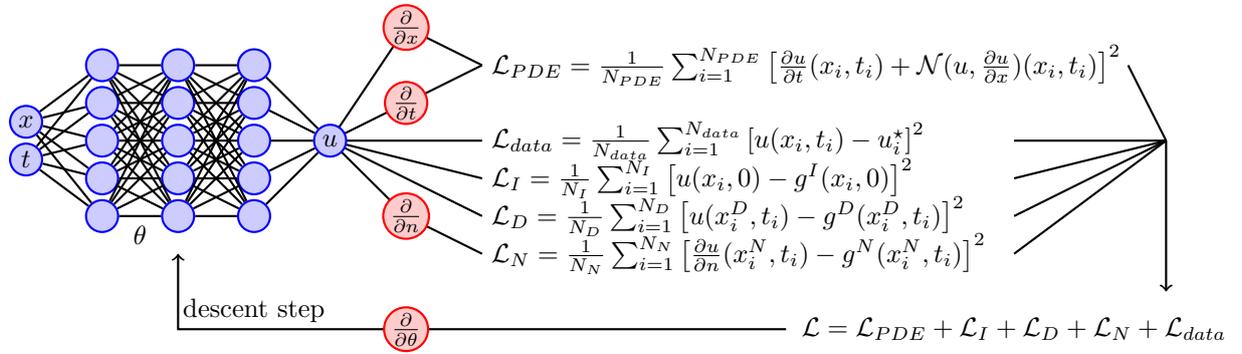

In this context, the input of the MLP network is composed of the space-time coordinates $(x,t)$ of a given point, while its output predicts the solution field $u$ at the same point.  The MLP network is controlled by a set of parameters $\theta$, which have to be calibrated to fulfill the governing equations and fit the observation data. This \emph{learning} phase consists in the minimization of a \emph{loss function} $\mathcal{L}$ that embeds all different criteria. An overview of the PINNs principle is depicted in Fig.~\ref{fig:pinns}, while the following subsections describe the different processing steps.

\begin{remark}
Additional physical parameters can be added as input of the network, allowing to construct a parametric surrogate, as will be shown in section~\ref{sec:param}.
\end{remark}

\begin{remark}
In the learning typology, PINNs can be considered as semi-supervised learning algorithms, since they rely partially on labelled data.
\end{remark}

\subsection{Multi-layer perceptron}

As explained above, the MLP network is the core of the PINNs methodology which predicts the solution of the problem. It is composed of a set of L layers: an input layer, $L-2$ intermediate hidden layers and an output layer. Each layer contains a set of $N_l$ neurons, each of them being connected to all the neurons of the previous layer. A neuron is a mathematical operator that applies a non-linear \emph{activation function} $\alpha$ to the weighted sum of its own inputs plus a bias factor. Thus, the value of the $i$th neuron of the layer $l$ writes:
\begin{equation}
f_i^l = \alpha(\sum_{j=1}^{N_{l-1}} w^l_{ij} f_j^{l-1} + b_i^l)
\end{equation}
where $w^l_{ij}$ is the weight of the connection with the neuron $j$ of the layer $l-1$ and $b_i^l$ is the bias. The values of the neurons in the first layer correspond to the MLP inputs.
As a consequence, the output of the MLP $u$ is obtained by propagating forward the input values $(x,t)$ through the network layers, resulting in a composition of the activation functions.
Such a network is characterized by its weights and biases, which constitute the set of parameters $\theta$ to be calibrated during the learning phase.
It has been shown that such networks have the ability to represent a large variety of functions~\cite{hornikMultilayerFeedforwardNetworks1989}. However, there is no guarantee that the learning procedure converges to a satisfactoy network~\cite{krishnapriyanCharacterizingPossibleFailure2021}.

\begin{remark}
    No activation is usually applied to the last layer, then the predicted solution is just a linear combination of the $N_{L-1}$ values obtained in the last but one layer. The use of the hidden layers can therefore be interpreted as the construction of $N_{L-1}$ suitable basis functions to represent the solution.
\end{remark}

\begin{remark}
    The ability of neural networks to represent complex functions depends on the number of neurons and their connections. Thus, the network parameters, i.e. connection weights and neuron biases, play the role of degrees of freedom for classical solvers. However, they are not located spatially, contrary to the case of nodal numerical schemes.
\end{remark}

\subsection{Physics-based loss}

To train the network, one seeks for the parameters $\theta$ that minimize a loss function $\mathcal{L}$, which should embed criteria reflecting both the fitting of observation data and the fulfillment of the governing equations with boundary and initial conditions. Thus, it is written as a sum of specific losses:
\begin{equation}
    \mathcal{L} = \mathcal{L}_{PDE} + \mathcal{L}_{I} + \mathcal{L}_{D} + \mathcal{L}_{N} + \mathcal{L}_{data}
\end{equation}
The four first loss terms are computed as the Mean Squared Error (MSE) associated to the PDE residuals, initial and boundary conditions, evaluated for a space-time sampling of the domain, whereas the last term is computed as the MSE associated to observation data:
\begin{align}
    &\mathcal{L}_{PDE}  = \frac{1}{N_{PDE}} \sum_{i=1}^{N_{PDE}} \left [ \frac{\partial{u}}{\partial{t}}(x_i,t_i) + \mathcal{N}(u,\frac{\partial{u}}{\partial{x}})(x_i,t_i) \right]^2 \\
    &\mathcal{L}_{D}  = \frac{1}{N_{D}} \sum_{i=1}^{N_{D}} \left [ u(x_i^D,t_i) - g^D(x_i^D,t_i) \right ]^2 \\
    &\mathcal{L}_{N} = \frac{1}{N_{N}} \sum_{i=1}^{N_{N}} \left [ \frac{\partial{u}}{\partial{n}}(x_i^N,t_i) - g^N(x_i^N,t_i) \right ]^2 \\
    &\mathcal{L}_{data} = \frac{1}{N_{data}} \sum_{i=1}^{N_{data}} \left [ u(x_i,t_i) - u^\star_i \right ]^2
\end{align}
Different types of sampling can be employed~\cite{wangExpertGuideTraining2023} but a random sampling based on a uniform or latin hypercube distribution usually performs well.

The evaluations of the derivatives of the solution (network output) with respect to the coordinates (network inputs), necessary to estimate the loss terms associated to the PDE residuals and the Neumann conditions, are carried out using automatic differentiation techniques, which are included in most AI software libraries.

\begin{remark}
    Several authors introduce a weighted sum of the losses to control the balance of the different terms~\cite{wangExpertGuideTraining2023}, but we found not mandatory to consider such additional parameters.
\end{remark}

\begin{remark}
    The use of exact derivatives is a major difference compared to classical PDE solvers, which approximate the differential operators to construct the numerical schemes.
\end{remark}

\begin{remark}
    This differentiation step necessitates the use of regular activation functions. Thus, some functions commonly employed in image processing (e.g. $Relu$) cannot be employed here, due to their lack of regularity.
\end{remark}

\begin{remark}
    The sampling plays a role similar to that of the mesh for classical solvers, since it defines where the PDE residuals are evaluated. However, sampling is not connected to the definition of the degrees of freedom. Note also that no approximation of the geometry is achieved during the sampling, contrary to the meshing step.
\end{remark}

\subsection{Learning procedure}

The minimization of the loss function $\mathcal{L}$ is the critical step of the learning task. This is a challenging optimization problem for the following reasons: i) the dimension of the variable $\theta$ can be extremely large ii) the non-linearity of the network representation yields non-convex, anisotropic and possibly multimodal loss functions.

Only gradient-based descent methods can solve such problems involving a large number of variables. Thus, automatic differentiation techniques are again used to compute the derivative of the loss function $\mathcal{L}$ with respect to the parameters $\theta$.

In machine learning, stochastic gradient methods are usually employed, due to the use of \emph{batch} procedures to split large datasets into smaller ones, which introduces randomness in the loss function evaluation. The most commonly employed method is ADAM algorithm, which is a first-order approach based on adaptive estimates of lower-order moments. However, several authors reported the low convergence rates obtained using ADAM algorithm to train PINNs, caused by the anisotropy of the loss function~\cite{krishnapriyanCharacterizingPossibleFailure2021,wangExpertGuideTraining2023}. As a consequence, higher-order descent methods, such as L-BFGS algorithm, are usually employed to accelerate the final convergence after an initial phase carried out using ADAM algorithm, which is more robust.

\begin{remark}
    L-BFGS algorithm and other quasi-Newton methods rely on the iterative construction of the Hessian matrix of the loss function and, therefore, do not comply with the use of batch procedures.
\end{remark}

\begin{remark}
    The specific form of the loss function, expressed as a sum of terms of different nature, yields ill-conditioned problems~\cite{krishnapriyanCharacterizingPossibleFailure2021}. Essentially, the minimization of all loss terms constitutes a multi-criterion optimization problem, which is difficult to solve by just summing the criteria.
\end{remark}

\subsection{Synthesis}

PINNs appear finally as a collocation method based on a neural network representation of the solution. Automatic differentiation plays a critical role, for the computation of the PDE residuals as well as the descent direction. The resulting procedure is particularly versatile and straightforward to implement, making the approach especially interesting for model experimentation. A second important characteristic is the ability of PINNs to handle both data and PDE models, filling the gap between simulation and experiments. However, PINNs do not rely so far on a safe theoretical basis to guarantee the results in terms of convergence and accuracy. Especially, the minimization of the loss function remains a difficult task, often problem dependent, despite the numerous algorithmic extensions proposed by the community~\cite{cuomoScientificMachineLearning2022,wangExpertGuideTraining2023}.

An interesting feature of PINNs is the use of an optimization formulation to solve the PDE system, yielding a very flexible approach that can be easily and efficiently extended to more complex problems. In the following sections, we show how to leverage this formulation to construct parametric surrogate models (see section~\ref{sec:param}), simulate multi-physic systems (see section~\ref{sec:coupling}) or assimilate data in the context of turbulence (see section~\ref{sec:assimilation}).

%%%%%%%%%%%%%%%%%%%%%%%%%%%%%%%%%%%%%%%%%%%%%%%%%%%%%%%%%%%%%%%%%%%%%%

\section{Parametric simulation}
\label{sec:param}

\subsection{Method}

The baseline PINNs methodology presented in the previous section can be easily extended to simulate a problem for a range of physical parameters, yielding a parametric surrogate. As example, we present in the following subsection the simulation of the flow in a differentially heated square cavity for a range of viscosity and conduction coefficients. Once the network is trained, any configuration can then be simulated for a negligible cost. This parametric approach is especially interesting for uncertainty quantification or design exploration.

The modification to carry out is straightforward since it just consists in introducing the physical parameters as additional inputs of the network and, consistently, extend the sampling to cover the domain. As illustration, we introduce a model parameter $\gamma$ affecting the spatial differential operator in the problem defined by Eq.~\ref{eq:pb}:
\begin{equation}
    \frac{\partial{u}}{\partial{t}} + \mathcal{N}_\gamma(u,\frac{\partial{u}}{\partial{x}}) = 0 \quad (x,t,\gamma) \in \Omega \times [0,T] \times [\gamma_m,\gamma_M]
\label{eq:pb_param}
\end{equation}
The procedure is then modified as follows: $\gamma$ is added as third input of the network and the evaluation of the loss terms is achieved for sampling points of extended coordinates $(x_i,t_i,\gamma_i)$. This modification is negligible in terms of implementation because the simulation is formulated as an optimization problem, whose only modification concerns the loss function evaluation. Nevertheless, the computational effort is obviously increased, due to the enlargement of the sampling.

\subsection{Application to heated cavity case}

We demonstrate the ability of PINNs to construct such parametric simulations for a natural convection problem with variations of two fluid coefficients: we consider a steady flow in a squared domain of size $2 \times 2$ subject to differentially heated lateral walls. The temperature of the left wall is $T_L = 1$, whereas the one of the right wall is $T_R = -1$, adiabaticity being imposed for top and bottom walls. The governing equations are the incompressible Navier-Stokes equations with boyancy and gravity, denoted as $\mathcal{N}_{NS}$:
\begin{align}
\label{eq:NS_boyancy_beg}
    & \frac{\partial u}{\partial x} + \frac{\partial v}{\partial y} = 0 \\
    & u \frac{\partial u}{\partial x} + v \frac{\partial u}{\partial y} + \frac{1}{\rho} \frac{\partial p}{\partial x} - \frac{\mu}{\rho} (\frac{\partial^2 u}{\partial x^2}+\frac{\partial^2 u}{\partial y^2}) = 0 \\
    & u \frac{\partial v}{\partial x} + v \frac{\partial v}{\partial y} + \frac{1}{\rho} \frac{\partial p}{\partial y} - \frac{\mu}{\rho} (\frac{\partial^2 v}{\partial x^2}+\frac{\partial^2 v}{\partial y^2}) - \beta g T -g = 0 \\
    & u \frac{\partial T}{\partial x} + v \frac{\partial T}{\partial y} - \frac{k^f}{\rho C_p} (\frac{\partial^2 T}{\partial x^2}+\frac{\partial^2 T}{\partial y^2}) = 0
\label{eq:NS_boyancy_end}
\end{align}
where the velocity is denoted as $(u,v)$, the pressure $p$ and the temperature $T$. The coefficients characterizing the fluid are the density $\rho = 1$, the viscosity $\mu$, the expansion coefficient $\beta=0.1$, the thermal conductivity $k^f$ and the heat capacity $C_p=1$. The objective is to simulate the flow for a range of viscosity $\mu \in [0.1,0.01]$ and conductivity $k^f \in [0.1,0.01]$.

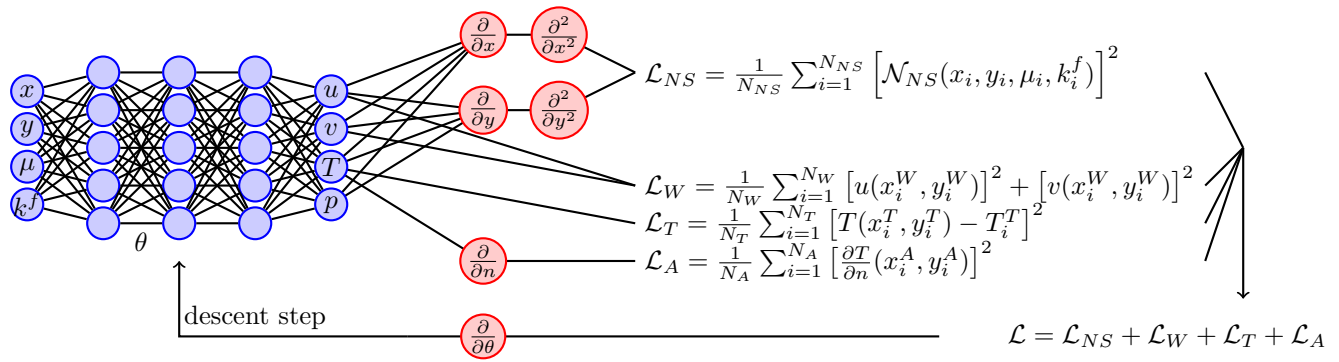
\begin{figure}[!ht]
\centering
\begin{tikzpicture}[x=1cm,y=0.5cm]

\node at (4,-4.5)  {descent step};
\draw[->,thick] (6,-5) -- (3,-5) -- (3,-3);

\draw[-,thick] (6,-5) -- (13,-5);
\node[mynode2] at (7,-5) {$\frac{\partial}{\partial{\theta}}$};

\draw[-,thick] (16.5,2) -- (17,0);
\draw[-,thick] (16.5,-1) -- (17,0);
\draw[-,thick] (16.5,-2) -- (17,0);
\draw[-,thick] (16.5,-3) -- (17,0);
\draw[->,thick] (17,0) -- (17,-4);
\node at (16,-5)  {$\mathcal{L} = \mathcal{L}_{NS} + \mathcal{L}_{W} + \mathcal{L}_{T} + \mathcal{L}_{A}$} ;

\draw[-,thick] (7,-3) -- (9,-3);
\draw[-,thick] (5,1.5) -- (9,-1);
\draw[-,thick] (5,0.5) -- (9,-1);
\draw[-,thick] (5,-0.5) -- (9,-2);

\node[right] at (9,-1)  {$\mathcal{L}_{W} = \frac{1}{N_{W}} \sum_{i=1}^{N_{W}} \left [ u(x_i^W,y_i^W) \right ]^2 + \left [ v(x_i^W,y_i^W) \right ]^2$};
\node[right] at (9,-2)  {$\mathcal{L}_{T} = \frac{1}{N_{T}} \sum_{i=1}^{N_{T}} \left [ T(x_i^T,y_i^T) - T^T_i \right ]^2$};
\node[right] at (9,-3)  {$\mathcal{L}_{A} = \frac{1}{N_{A}} \sum_{i=1}^{N_{A}} \left [ \frac{\partial T}{\partial n}(x_i^A,y_i^A) \right ]^2$};

\draw[-,thick] (8,3) -- (9,2);
\draw[-,thick] (8,1) -- (9,2);
\node[right] at (9,2)  {$\mathcal{L}_{NS} = \frac{1}{N_{NS}} \sum_{i=1}^{N_{NS}} \left [
\mathcal{N}_{NS}(x_i,y_i,\mu_i,k^f_i)
 \right]^2$};

\draw[-,thick] (5,1.5) -- (7,3);
\draw[-,thick] (5,0.5) -- (7,3);
\draw[-,thick] (5,-0.5) -- (7,3);
\draw[-,thick] (5,-1.5) -- (7,3);
\draw[-,thick] (5,1.5) -- (7,1);
\draw[-,thick] (5,0.5) -- (7,1);
\draw[-,thick] (5,-0.5) -- (7,1);
\draw[-,thick] (5,-1.5) -- (7,1);

\draw[-,thick] (7,3) -- (8,3);
\draw[-,thick] (7,1) -- (8,1);

\draw[-,thick] (5,-0.5) -- (7,-3);
\node[mynode2] at (7,3)  {$\frac{\partial}{\partial x}$};
\node[mynode2] at (7,1)  {$\frac{\partial}{\partial y}$};
\node[mynode2] at (8,3)  {$\frac{\partial^2}{\partial x^2}$};
\node[mynode2] at (8,1)  {$\frac{\partial^2}{\partial y^2}$};
\node[mynode2] at (7,-3)  {$\frac{\partial}{\partial{n}}$};

\node at (2.5,-2.5)  {$\theta$};
  \readlist\Nnod{4,5,5,5,4} % number of nodes per layer
  \foreachitem \N \in \Nnod{ % loop over layers
    \foreach \i [evaluate={\x=\Ncnt; \y=\N/2-\i+0.5; \prev=int(\Ncnt-1);}] in {1,...,\N}{ % loop over nodes
      \node[mynode] (N\Ncnt-\i) at (\x,\y) {};
      \ifnum\Ncnt>1 % connect to previous layer
        \foreach \j in {1,...,\Nnod[\prev]}{ % loop over nodes in previous layer
          \draw[-,thick] (N\prev-\j) -- (N\Ncnt-\i); % connect arrows directly
        }
      \fi % else: nothing to connect first layer
    }
  }
\node at (1,1.5)  {$x$};
\node at (1,0.5)  {$y$};
\node at (1,-0.5)  {$\mu$};
\node at (1,-1.5)  {$k^f$};
\node at (5,1.5)  {$u$};
\node at (5,0.5)  {$v$};
\node at (5,-0.5)  {$T$};
\node at (5,-1.5)  {$p$};

\end{tikzpicture}
\caption{Parametric simulation for the heated cavity case}
\label{fig:pinns_param}
\end{figure}

\begin{figure}[!ht]
\centering
\includegraphics[width=0.4\linewidth]{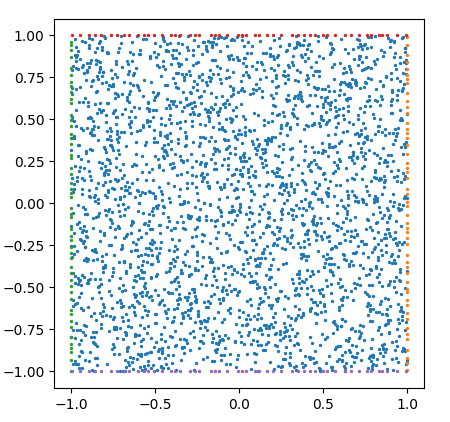}
\caption{Sampling for the cavity case.}
\label{fig:cavity_sampling}
\end{figure}

The network considered has four inputs $(x,y,\mu,k^f)$ and four outputs $(u,v,T,p)$, as depicted in Fig.~\ref{fig:pinns_param}. It includes 3 hidden layers of 50 neurons, with hyperbolic tangent activation functions $\alpha = tanh$. The loss function is the sum of the MSE terms related to the residuals of the governing equations $\mathcal{L}_{NS}$, the no-slip conditions on the walls $\mathcal{L}_W$, the adiabatic conditions on the top and botton walls $\mathcal{L}_A$ and the imposed temperature on the left and right walls $\mathcal{L}_T$. Note that no observation data is provided for this case. Regarding the sampling, a latin hypercube distribution of size $2500$ is adopted for the spatial domain, as shown in Fig.~\ref{fig:cavity_sampling}. Additionnally, the range of variation of the coefficients $\nu$ and $k^f$ is discretized with 4 equidistributed values. Thus, the complete sampling based on a tensorial product counts $2500 \times 4 \times 4$ points. The minimization of the loss function is achieved using the ADAM algorithm (3,000 epochs, learning rate $10^{-3}$) followed by the L-BFGS algorithm (7,000 epochs).

\begin{figure}[!ht]
\centering
\includegraphics[width=0.33\linewidth]{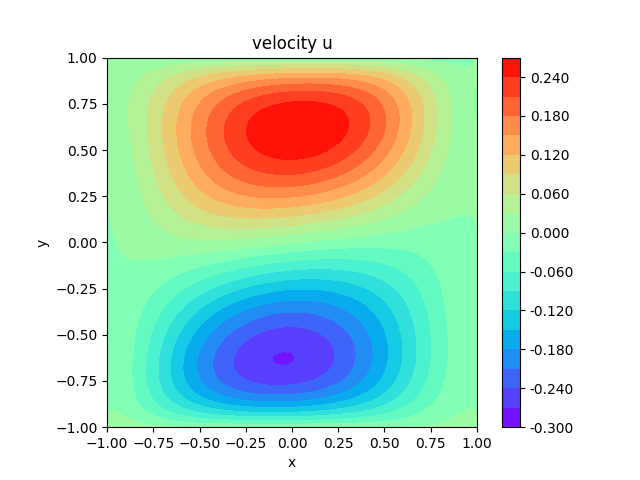}\hfill
\includegraphics[width=0.33\linewidth]{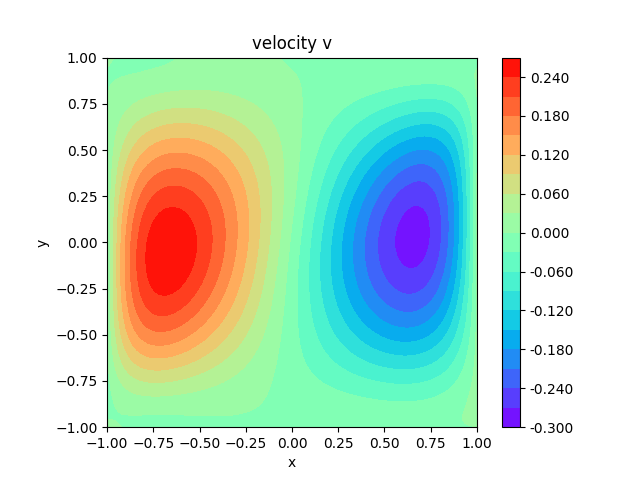}\hfill
\includegraphics[width=0.33\linewidth]{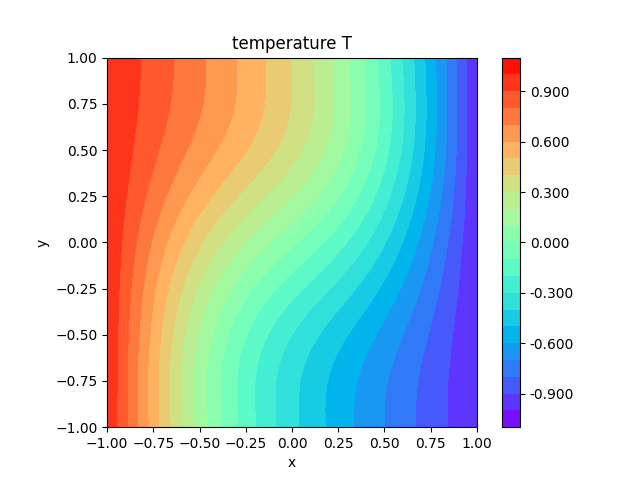}\\
\includegraphics[width=0.33\linewidth]{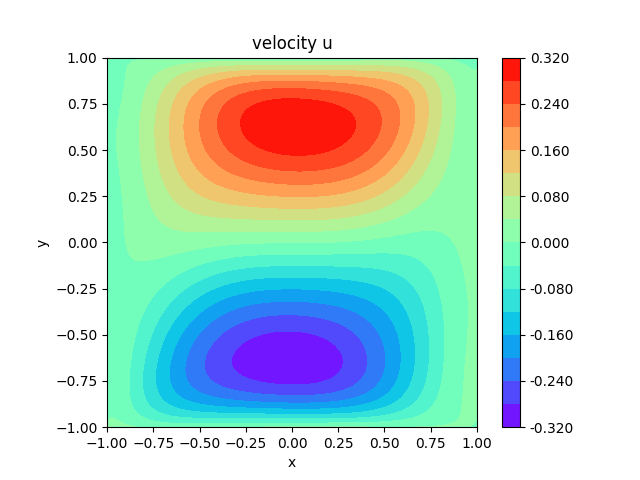}\hfill
\includegraphics[width=0.33\linewidth]{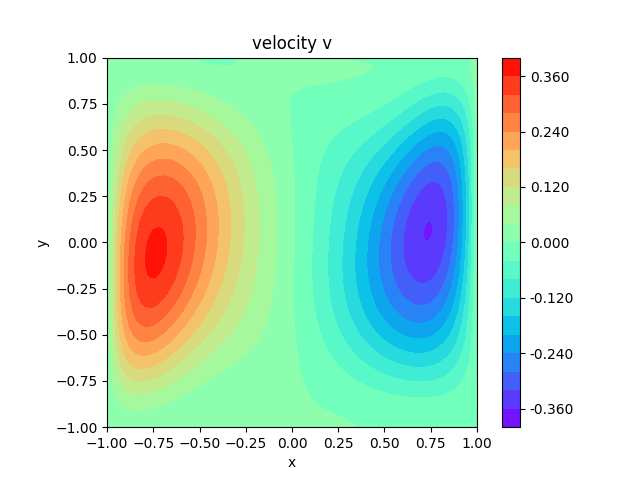}\hfill
\includegraphics[width=0.33\linewidth]{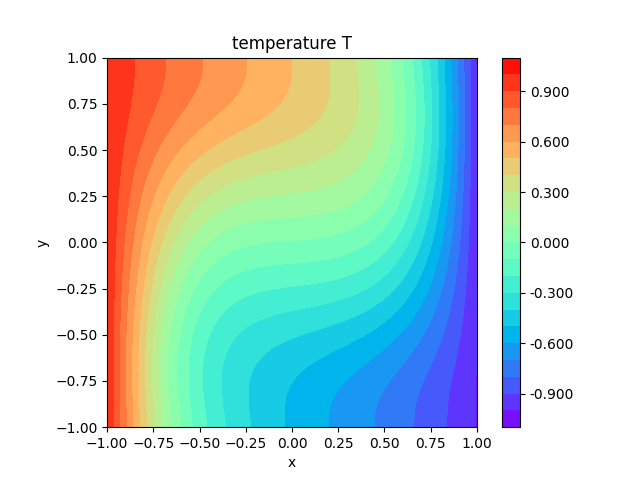}\\
\includegraphics[width=0.33\linewidth]{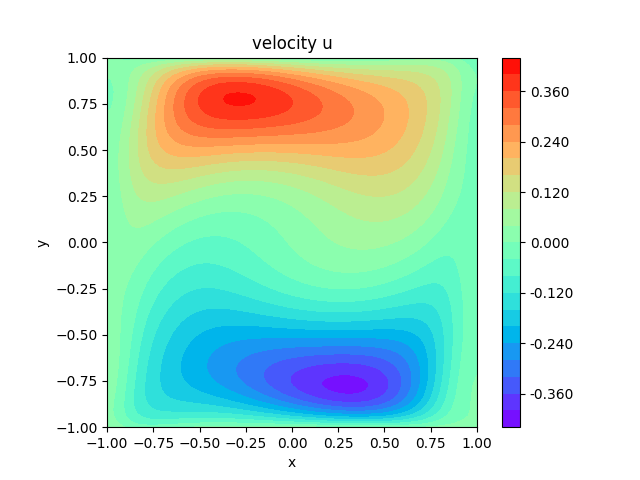}\hfill
\includegraphics[width=0.33\linewidth]{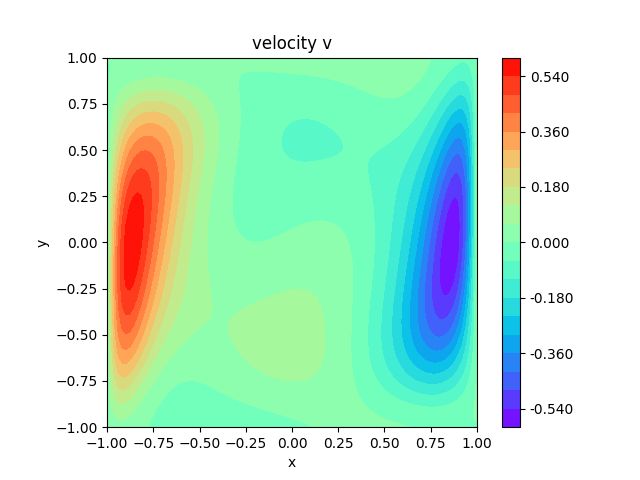}\hfill
\includegraphics[width=0.33\linewidth]{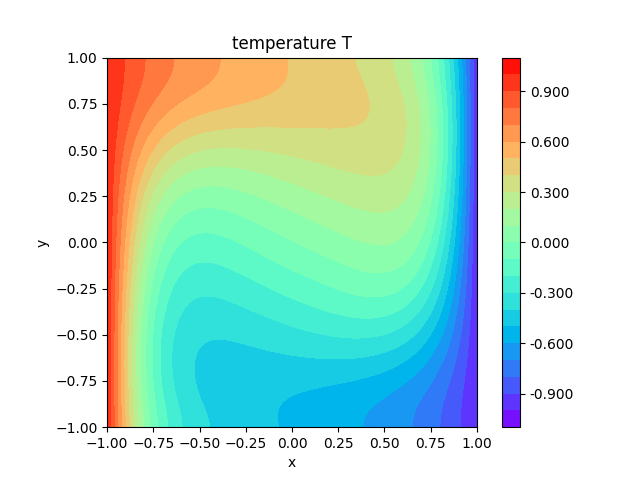}
\caption{Predictions for different parameter values (parametric simulation): $\mu=0.1$ and $k^f=0.1$ (top), $\mu=0.05$ and $k^f=0.05$ (middle), $\mu=0.01$ and $k^f=0.01$ (bottom)}
\label{fig:cavity_sol}
\end{figure}

\begin{figure}[!ht]
\centering
\includegraphics[width=0.33\linewidth]{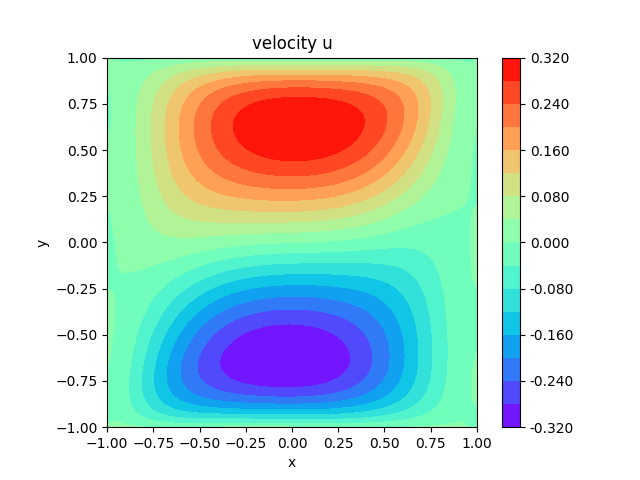}\hfill
\includegraphics[width=0.33\linewidth]{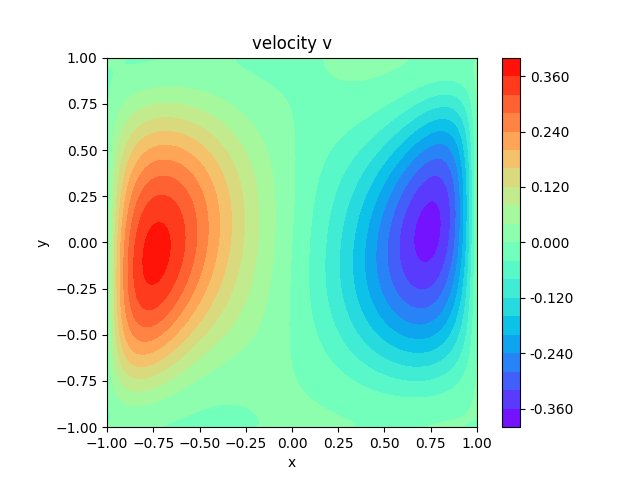}\hfill
\includegraphics[width=0.33\linewidth]{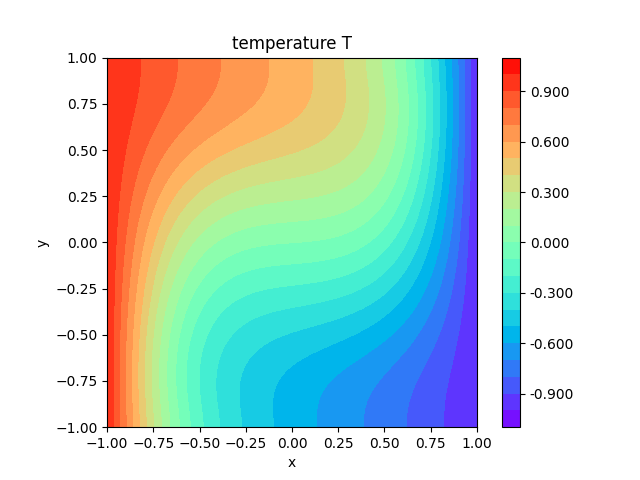}
\caption{Prediction for $\mu=0.05$ and $k^f=0.05$ (single-parameter simulation )}
\label{fig:cavity_test}
\end{figure}

The solution fields predicted for some parameter values can be seen in Fig.~\ref{fig:cavity_sol}. For the case $(\nu,k^f)=(0.05,0.05)$, a separate network has been trained independently to assess the efficiency of the parametric model. As shown in Fig.~\ref{fig:cavity_test}, a satisfactory agreement is observed between the fields obtained from the parametric model and the single-parameter model.

\clearpage

%%%%%%%%%%%%%%%%%%%%%%%%%%%%%%%%%%%%%%%%%%%%%%%%%%%%%%%%%%%%%%%%%%

\section{Multidisciplinary couplings}
\label{sec:coupling}

\subsection{Method}

The coupling of different disciplines, encountered for instance in aero-thermal or aero-structural problems, is generally a difficult task because of the presence of phenomena characterized by different spatial and temporal scales and by different mathematical natures (e.g. hyperbolic vs elliptic). As a consequence, each discipline is simulated by using a specific mesh, time step and numerical method, yielding a tedious coupling procedure at the interface.

In contrast, the PINNs formulation allows a straightforward implementation for coupled systems: each discipline is predicted by its own network, thus permitting to adjust the network complexity to the concerned physics, whereas a global loss function gathers the contributions of the different PDEs, boundary conditions, possible data, as well as coupling conditions. By minimizing this global loss, one solves simultaneously the different physics and their coupling. Of particular interest is the avoidance of any fixed-point algorithm to ensure the convergence of the coupling conditions, which would necessitate to implement an additional iterative loop. Note also that all coupling conditions are accounted simultaneously and it is not necessary to devise a specific strategy (e.g. Dirichlet-Neumann method) to establish the coupling.

In the following section, a typical implementation is detailed in the contexte of a conjugate heat transfer problem.

\subsection{Application to conjugate heat transfer}

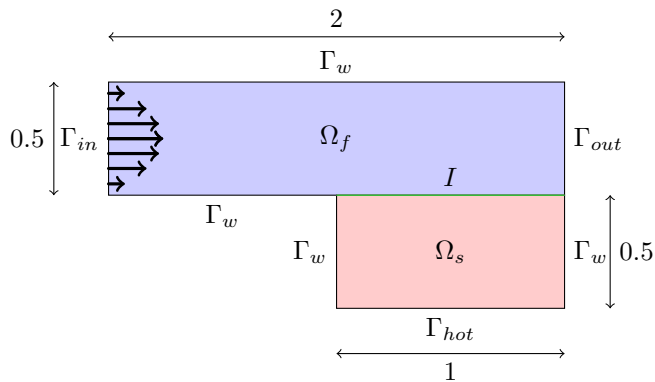
\begin{figure}[ht]
\begin{center}
    \begin{tikzpicture}[scale=0.6]

    % Fluid domain
    \draw [fill=blue!20, even odd rule] (0,0) rectangle (10,2.5)
        node[left]  at (0,1.25)  {$\Gamma_{in}$}
        node[above] at (5,2.5)   {$\Gamma_w$}
        node[right] at (10,1.25) {$\Gamma_{out}$}
        node[below] at (2.5,0)   {$\Gamma_w$}
        node[midway] {$\Omega_f$};

    % Fluid lengths
    \draw[-latex, <->] (0,3.5) -- (10,3.5) node[midway,above] {$2$};
    \draw[-latex, <->] (-1.2,0) -- (-1.2,2.5) node[midway,left] {$0.5$};

    % Solid domain
    \draw [fill=red!20, even odd rule] (5,0) rectangle (10,-2.5)
        node[left]  at (5,-1.25)  {$\Gamma_w$}
        node[right] at (10,-1.25) {$\Gamma_w$}
        node[below] at (7.5,-2.5) {$\Gamma_{hot}$}
        node[above] at (7.5,0)    {$I$}
        node[midway] {$\Omega_s$};

    % Solid lengths
    \draw[-latex, <->] (11,0) -- (11, -2.5) node[midway,right] {$0.5$};
    \draw[-latex, <->] (5,-3.5) -- (10,-3.5) node[midway,below] {$1$};

    \draw[white] (5,0)  -- (10,0);
    \draw[green!70!black] (5,0)  -- (10,0);

    \def\Ry{1} % big pipe vertical radius right
    \def\v{1.2}   % velocity magnitude

    \foreach \y [evaluate={ \vy=\v*(1-0.7*\y^2);}] in {0,0.33,0.66,1.0}{
    \draw[->,very thick,line cap=round] (0, \y*\Ry + 1.25) --++ (\vy,0);
    \draw[->,very thick,line cap=round] (0, -\y*\Ry+ 1.25) --++ (\vy,0);
    }

    \end{tikzpicture}
\end{center}
\caption{Conjugate heat transfer problem}
\label{fig:conjugate}
\end{figure}

The two-dimensional computational domain is composed of a fluid part $\Omega_f$ and a solid part $\Omega_s$, delimited by a coupling interface $I$, as depicted in Fig.~\ref{fig:conjugate}. Regarding the fluid, we consider a channel with imposed velocity (parabolic distribution with maximum $u_{in}=1$) and temperature $T_{in}=0.2$ at inlet $\Gamma_{in}$, whereas a flow rate conservation is imposed at outlet $\Gamma_{out}$. A no-slip condition is prescribed at the walls $\Gamma_{w}$ and $I$. Adiabaticity conditions are imposed at the walls $\Gamma_{w}$.

The flow is governed by incompressible Navier-Stokes equations with thermal transport~\ref{eq:NS_boyancy_beg}-\ref{eq:NS_boyancy_end}, but neglecting the gravity effects. The fluid coefficients are $\rho=1$, $\mu=0.01$, $k^f=0.025$ and $C_p=1$. Regarding the solid, the temperature $T_{hot}=1$ is fixed at the bottom boundary $\Gamma_{hot}$, whereas adiabaticity is imposed on side walls $\Gamma_{w}$. The temperature field in the solid is governed by the heat transfer equation, denoted as $\mathcal{N}_{HT}$:
\begin{equation}
    k^s (\frac{\partial^2 T}{\partial x^2}+\frac{\partial^2 T}{\partial y^2}) = 0
\end{equation}
where $k^s=1$ is the conductivity coefficient in the solid. At the interface $I$ between fluid and solid domains, a two-condition coupling is adopted, expressing the continuity of the temperature and thermal flux:
\begin{align}
    \label{eq:coupling_beg}
    T |_f &= T |_s \\
    k^f \left. \frac{\partial T}{\partial n} \right |_f &=  k^s \left. \frac{\partial T}{\partial n} \right |_s
    \label{eq:coupling_end}
\end{align}

\begin{figure}[!ht]

\begin{tikzpicture}[x=1cm,y=0.5cm]

\draw[->,thick] (16,-7.5) -- (16,-8.5);
\draw[->,thick] (15,-9) -- (15.5,-9);

\draw[-,thick] (4,-0.5) -- (6,-8) -- (7,-8);
\draw[-,thick] (6,-5) -- (7,-10);

\node[right] at (7,-8)  {$\mathcal{L}_{c_1} = \frac{1}{N_{c}} \sum_{i=1}^{N_{c}} \left [ \left. T(x_i^{c},y_i^{c}) \right |_f - \left. T(x_i^{c},y_i^{c}) \right |_s \right ]^2$};
\node[right] at (7,-10)  {$\mathcal{L}_{c_2} = \frac{1}{N_{c}} \sum_{i=1}^{N_{c}} \left [ \left. k^f \frac{\partial T}{\partial n} \right |_f (x_i^{c},y_i^{c}) - \left. k^s \frac{\partial T}{\partial n} \right |_s (x_i^{c},y_i^{c}) \right ]^2$};

\node at (17,-9)  {$\mathcal{L} = \mathcal{L}_{f} + \mathcal{L}_{s} + \mathcal{L}_{c}$} ;

\node at (4,4.5)  {descent step};
\draw[->,thick] (5,5) -- (2.5,5) -- (2.5,3);

\draw[-,thick] (5,5) -- (18,5) -- (18,-8);
\node[mynode2] at (6,5) {$\frac{\partial}{\partial{\theta^f}}$};

\draw[-,thick] (15.5,2) -- (16,0);
\draw[-,thick] (15.5,-1) -- (16,0);
\draw[-,thick] (15.5,-2) -- (16,0);
\draw[-,thick] (15.5,-3) -- (16,0);
\draw[-,thick] (15.5,-4) -- (16,0);
\draw[-,thick] (15.5,-5) -- (16,0);
\draw[->,thick] (16,0) -- (16,-6);
\node at (14,-6.8)  {$\mathcal{L}_f = \mathcal{L}_{NS} + \mathcal{L}_{W} + \mathcal{L}_{u_{in}} + \mathcal{L}_{q} + \mathcal{L}_{T_{in}} + \mathcal{L}_{A}$} ;

\draw[-,thick] (6,-5) -- (8,-5);
\draw[-,thick] (4,1.5) -- (8,-1);
\draw[-,thick] (4,0.5) -- (8,-1);
\draw[-,thick] (4,1.5) -- (8,-2);
\draw[-,thick] (4,1.5) -- (8,-3);
\draw[-,thick] (4,-0.5) -- (8,-4);

\node[right] at (8,-1)  {$\mathcal{L}_{W} = \frac{1}{N_{W}} \sum_{i=1}^{N_{W}} \left [ u(x_i^W,y_i^W) \right ]^2 + \left [ v(x_i^W,y_i^W) \right ]^2$};
\node[right] at (8,-2)  {$\mathcal{L}_{u_{in}} = \frac{1}{N_{in}} \sum_{i=1}^{N_{in}} \left [ u(x_i^{in},y_i^{in}) - u_{in} \right ]^2$};
\node[right] at (8,-3)  {$\mathcal{L}_{q} = \left [ \frac{1}{N_{out}} \sum_{i=1}^{N_{out}} u(x_i^{out},y_i^{out}) - q_{out}/\Gamma_{out} \right ]^2$};
\node[right] at (8,-4)  {$\mathcal{L}_{T_{in}} = \frac{1}{N_{in}} \sum_{i=1}^{N_{in}} \left [ T(x_i^{in},y_i^{in}) - T_{in} \right ]^2$};
\node[right] at (8,-5)  {$\mathcal{L}_{A} = \frac{1}{N_{A}} \sum_{i=1}^{N_{A}} \left [ \frac{\partial T}{\partial n}(x_i^A,y_i^A) \right ]^2$};

\draw[-,thick] (7,3) -- (8,2);
\draw[-,thick] (7,1) -- (8,2);
\node[right] at (8,2)  {$\mathcal{L}_{NS} = \frac{1}{N_{NS}} \sum_{i=1}^{N_{NS}} \left [
\mathcal{N}_{NS}(x_i,y_i)
 \right]^2$};

\draw[-,thick] (4,1.5) -- (6,3);
\draw[-,thick] (4,0.5) -- (6,3);
\draw[-,thick] (4,-0.5) -- (6,3);
\draw[-,thick] (4,-1.5) -- (6,3);
\draw[-,thick] (4,1.5) -- (6,1);
\draw[-,thick] (4,0.5) -- (6,1);
\draw[-,thick] (4,-0.5) -- (6,1);
\draw[-,thick] (4,-1.5) -- (6,1);

\draw[-,thick] (6,3) -- (7,3);
\draw[-,thick] (6,1) -- (7,1);

\draw[-,thick] (4,-0.5) -- (6,-5);
\node[mynode2] at (6,3)  {$\frac{\partial}{\partial x}$};
\node[mynode2] at (6,1)  {$\frac{\partial}{\partial y}$};
\node[mynode2] at (7,3)  {$\frac{\partial^2}{\partial x^2}$};
\node[mynode2] at (7,1)  {$\frac{\partial^2}{\partial y^2}$};
\node[mynode2] at (6,-5)  {$\frac{\partial}{\partial{n}}$};

\node at (2.5,2.5)  {$\theta^f$};
  \readlist\Nnod{2,5,5,4} % number of nodes per layer
  \foreachitem \N \in \Nnod{ % loop over layers
    \foreach \i [evaluate={\x=\Ncnt; \y=\N/2-\i+0.5; \prev=int(\Ncnt-1);}] in {1,...,\N}{ % loop over nodes
      \node[mynode] (N\Ncnt-\i) at (\x,\y) {};
      \ifnum\Ncnt>1 % connect to previous layer
        \foreach \j in {1,...,\Nnod[\prev]}{ % loop over nodes in previous layer
          \draw[-,thick] (N\prev-\j) -- (N\Ncnt-\i); % connect arrows directly
        }
      \fi % else: nothing to connect first layer
    }
  }
\node at (1,0.5)  {$x$};
\node at (1,-0.5)  {$y$};
\node at (4,1.5)  {$u$};
\node at (4,0.5)  {$v$};
\node at (4,-0.5)  {$T$};
\node at (4,-1.5)  {$p$};

\fill[blue,very nearly transparent] (1.5,-7.2) -- (17.5,-7.2) -- (17.5,4) -- (1.5,4) -- cycle;

\node[blue] at (3,-6)  {fluid model};
\node at (3,-9)  {coupling};

\end{tikzpicture}

%%%
\vspace*{-2cm}
\hspace*{-0.5cm}
%%%

\begin{tikzpicture}[x=1cm,y=0.5cm]

    \draw[->,thick] (16,5.5) -- (16,6.5);

    \draw[-,thick] (4,0) -- (6,8) -- (7,8);
    \draw[-,thick] (6,3) -- (7,6);

    \node at (15,4.8)  {$\mathcal{L}_s = \mathcal{L}_{HT} + \mathcal{L}_{hot} + \mathcal{L}_{A}$} ;

    \node at (4,-4.5)  {descent step};
    \draw[->,thick] (6,-5) -- (2.5,-5) -- (2.5,-2);

    \draw[-,thick] (6,-5) -- (18,-5) -- (18,6);
    \node[mynode2] at (6,-5) {$\frac{\partial}{\partial{\theta^s}}$};

    \draw[-,thick] (15.5,2) -- (16,0);
    \draw[-,thick] (15.5,3) -- (16,0);
    \draw[-,thick] (15.5,-2) -- (16,0);
    \draw[->,thick] (16,0) -- (16,4);

    \draw[-,thick] (6,3) -- (8,3);
    \draw[-,thick] (4,0) -- (8,2);

    \node[right] at (8,2)  {$\mathcal{L}_{hot} = \frac{1}{N_{hot}} \sum_{i=1}^{N_{hot}} \left [ T(x_i^{hot},y_i^{hot}) - T_{hot} \right ]^2$};
    \node[right] at (8,3)  {$\mathcal{L}_{A} = \frac{1}{N_{A}} \sum_{i=1}^{N_{A}} \left [ \frac{\partial T}{\partial n}(x_i^A,y_i^A) \right ]^2$};

    \draw[-,thick] (7,-3) -- (8,-2);
    \draw[-,thick] (7,-1) -- (8,-2);
    \node[right] at (8,-2)  {$\mathcal{L}_{HT} = \frac{1}{N_{HT}} \sum_{i=1}^{N_{HT}} \left [
    \mathcal{N}_{HT}(x_i,y_i)
     \right]^2$};

    \draw[-,thick] (4,0) -- (6,3);
    \draw[-,thick] (4,0) -- (6,-3);
    \draw[-,thick] (4,0) -- (6,-1);

    \draw[-,thick] (6,-3) -- (7,-3);
    \draw[-,thick] (6,-1) -- (7,-1);

    \node[mynode2] at (6,-3)  {$\frac{\partial}{\partial x}$};
    \node[mynode2] at (6,-1)  {$\frac{\partial}{\partial y}$};
    \node[mynode2] at (7,-3)  {$\frac{\partial^2}{\partial x^2}$};
    \node[mynode2] at (7,-1)  {$\frac{\partial^2}{\partial y^2}$};
    \node[mynode2] at (6,3)  {$\frac{\partial}{\partial{n}}$};

    \node at (2.5,-1.5)  {$\theta^s$};
      \readlist\Nnod{2,3,3,1} % number of nodes per layer
      \foreachitem \N \in \Nnod{ % loop over layers
        \foreach \i [evaluate={\x=\Ncnt; \y=\N/2-\i+0.5; \prev=int(\Ncnt-1);}] in {1,...,\N}{ % loop over nodes
          \node[mynode] (N\Ncnt-\i) at (\x,\y) {};
          \ifnum\Ncnt>1 % connect to previous layer
            \foreach \j in {1,...,\Nnod[\prev]}{ % loop over nodes in previous layer
              \draw[-,thick] (N\prev-\j) -- (N\Ncnt-\i); % connect arrows directly
            }
          \fi % else: nothing to connect first layer
        }
      }
    \node at (1,0.5)  {$x$};
    \node at (1,-0.5)  {$y$};
    \node at (4,0)  {$T$};

    \fill[red,very nearly transparent] (1.5,-4) -- (17.5,-4) -- (17.5,5.2) -- (1.5,5.2) -- cycle;

    \node[red] at (3,4)  {solid model};

\end{tikzpicture}
\caption{Coupled simulation for the conjugate heat tranfer problem}
\label{fig:coupling}
\end{figure}
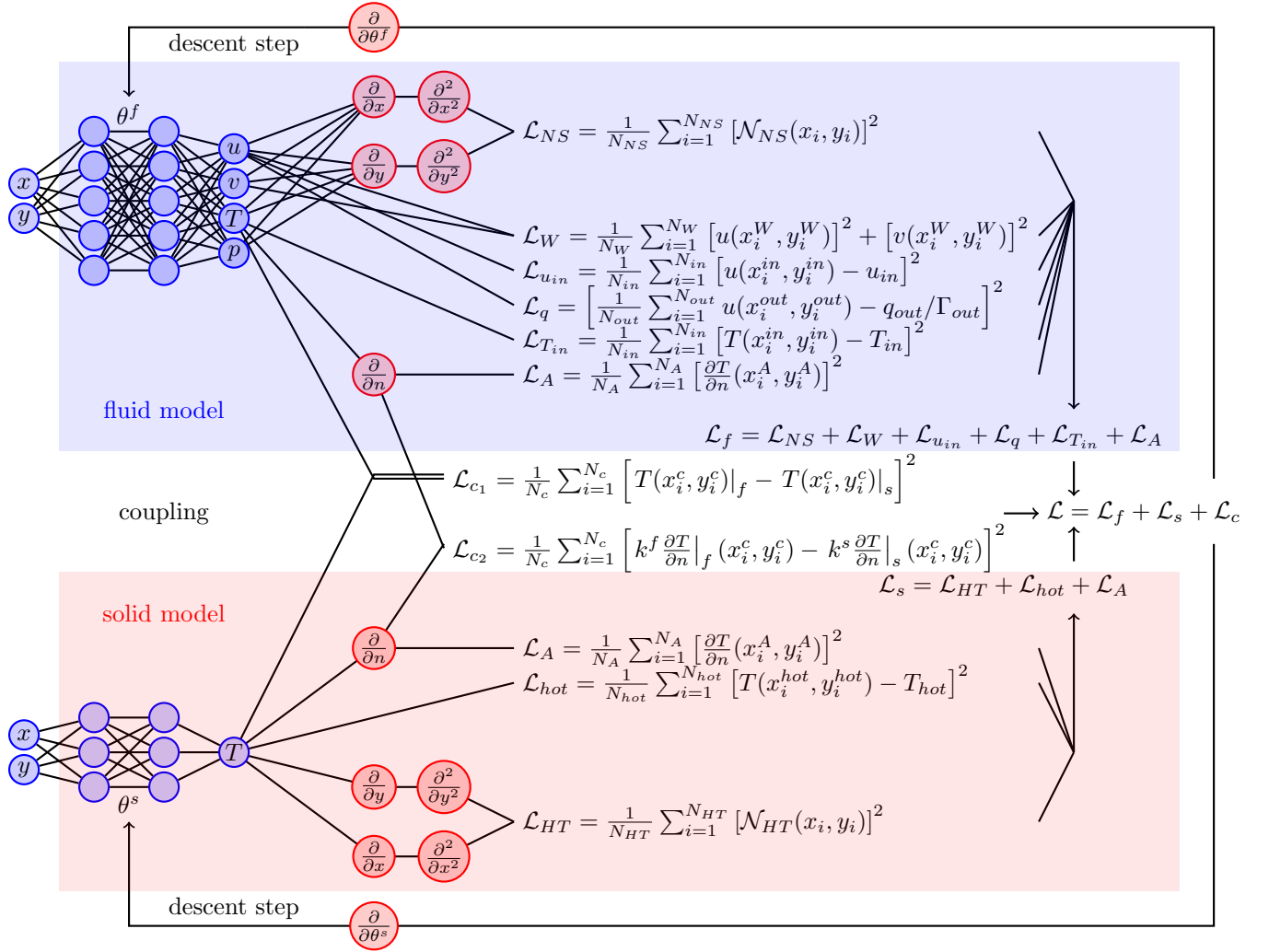

The simulation relies on two networks, with two inputs $(x,y)$, and respectively four outputs $(u,v,T,p)$ for the fluid one and a single output $T$ for the solid one, as shown in Fig.~\ref{fig:coupling}. The networks count respectively three hidden layers of 50 neurons and three hidden layers of 20 neurons, both of them employing the hyperbolic tangent activation function.

The domains are sampled using respectively 1600 and 225 points. For the fluid part, the training requires the computation of the loss term related to the Navier-Stokes equations $\mathcal{L}_{NS}$, the no-slip condition at the wall $\mathcal{L}_{W}$, the inlet condition for the velocity $\mathcal{L}_{u_{in}}$, the outlet condition for the flow rate $\mathcal{L}_{q}$, the adiabaticity condition $\mathcal{L}_{A}$ and the imposed temperature at inlet $\mathcal{L}_{T_{in}}$.
For the solid part, one has to evaluate the loss terms related to the heat transfer equation $\mathcal{L}_{HT}$, the adiabaticity condition $\mathcal{L}_{A}$ and the imposed temperature at the hot boundary $\mathcal{L}_{hot}$.
Moreover, the coupling conditions are implemented as two additional loss terms $\mathcal{L}_{C_1}$ and $\mathcal{L}_{C_2}$ expressing the MSE of the conditions~\ref{eq:coupling_beg}-\ref{eq:coupling_end}. These terms use the predictions of both networks. The global loss is finally evaluated by sommation of all the terms and differentiated with respect to the parameters of each network $\theta^f$ and $\theta^s$. The update of the two networks can then be achieved independently. As clearly shown in Fig.~\ref{fig:coupling}, a single iterative loop is performed to solve simultaneously the physics related to the fluid and solid, as well as the coupling conditions.
The training is performed using the ADAM algorithm for 50,000 epochs (learning rate $0.001$) followed by the BFGS algorithm for 2000 epochs. Note that other optimization strategies can be considered to update the networks~\cite{coulaud:hal-04225990}.

\begin{figure}[!ht]
\begin{minipage}{0.5\linewidth}
\centering
\includegraphics[width=0.9\linewidth]{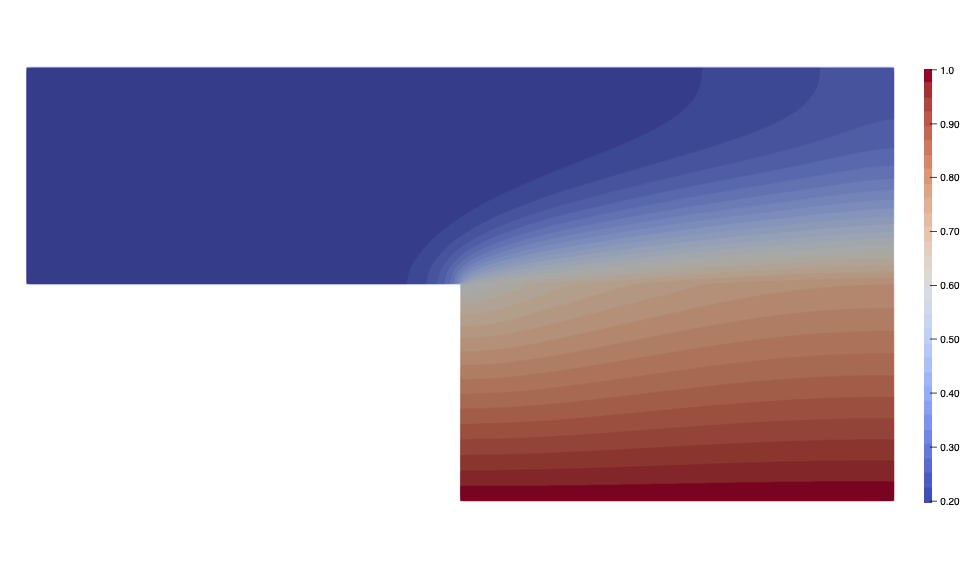}
\end{minipage}\hfill
\begin{minipage}{0.5\linewidth}
\centering
\includegraphics[width=0.9\linewidth]{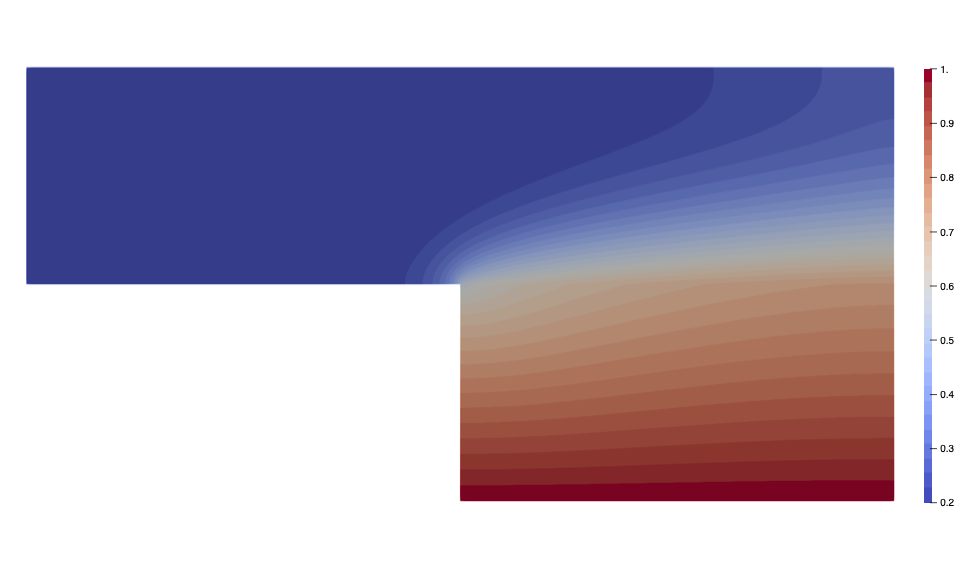}
\end{minipage}
\caption{Temperature predicted: left PINNs, right Finite-Element method}
\label{fig:temperature}
\end{figure}

\begin{figure}[!ht]
\begin{minipage}{0.5\linewidth}
\centering
\includegraphics[width=0.9\linewidth]{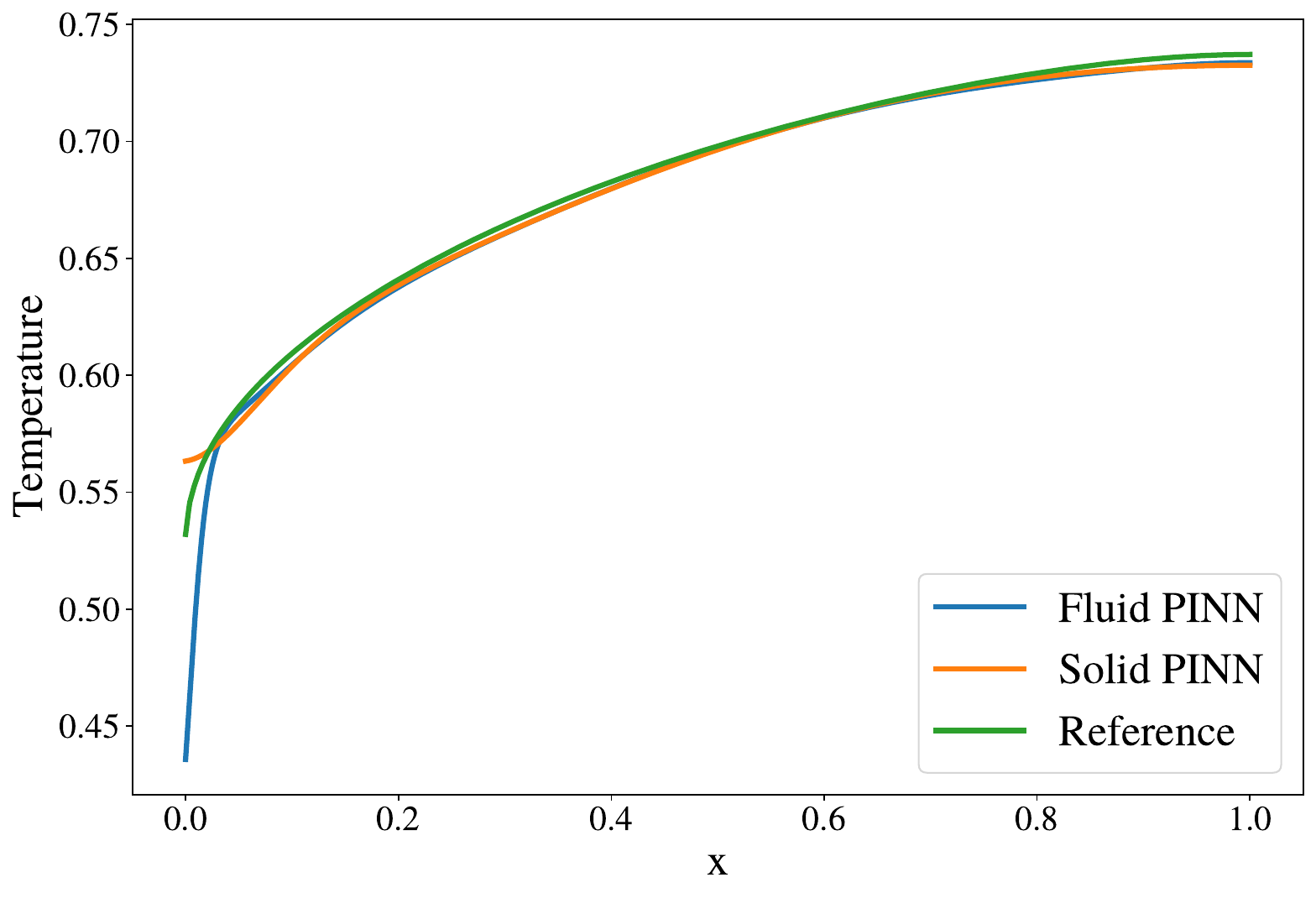}
\end{minipage}\hfill
\begin{minipage}{0.5\linewidth}
\centering
\includegraphics[width=0.9\linewidth]{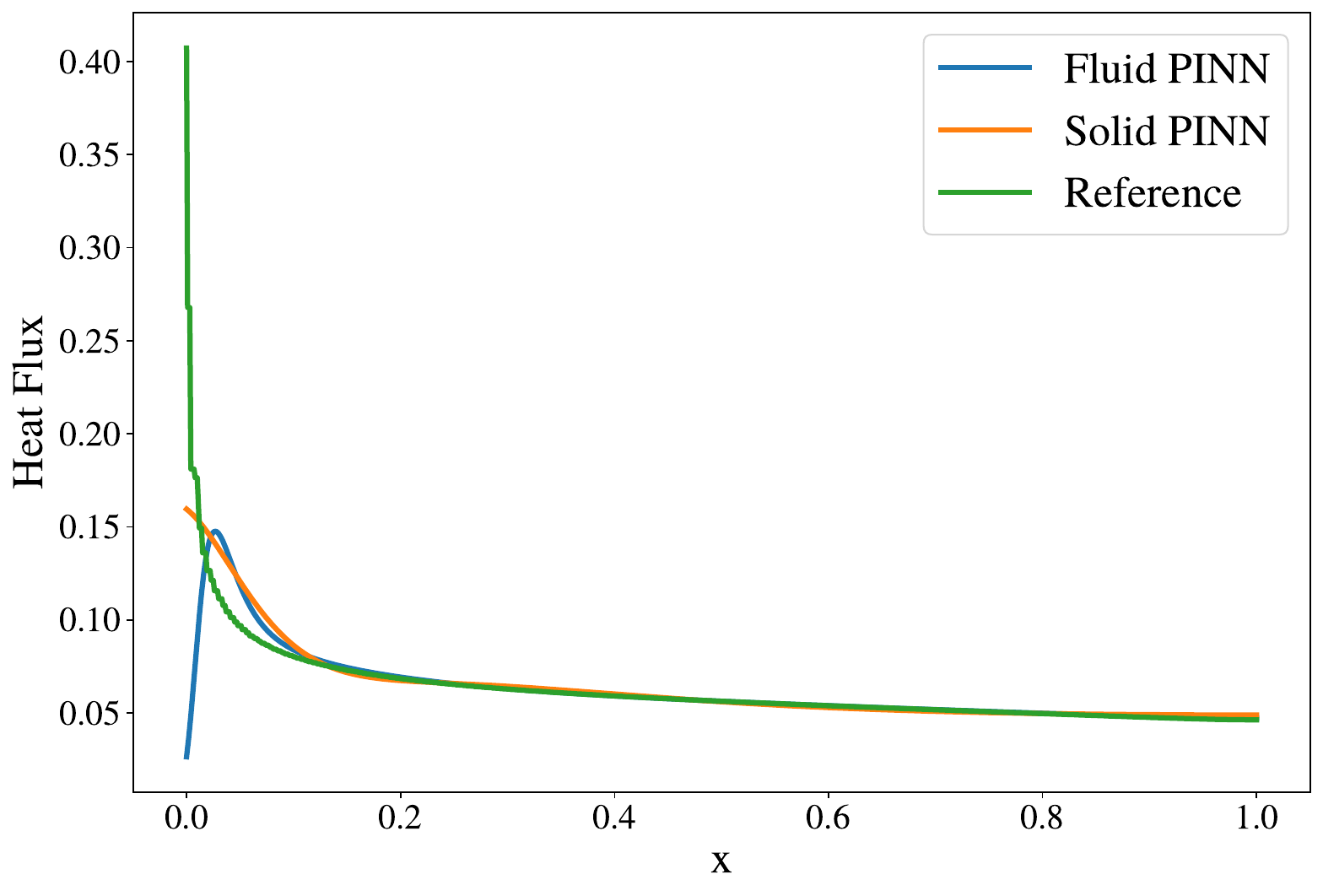}
\end{minipage}
\caption{Temperature (left) and heat flux (right) at the interface}
\label{fig:joint_interface_flux}
\end{figure}

The temperature field obtained using the PINNs approach is compared to the one resulting from a Finite-Element computation based on FreeFEM software\footnote{https://freefem.org/}, as shown in Fig~\ref{fig:temperature}. The agreement is satisfactory, however some discrepancies can be observed at the coupling interface. The temperature and heat flux at the interface are plotted in Fig.~\ref{fig:joint_interface_flux} and compared to the FreeFEM reference. As seen, the two coupling conditions are correct along the interface, except at the origin where the network is not able to capture the flux discontinuity, due to the regularity of the activation function employed.
Finally, some metrics are computed on a fine cartesian grid of size $100 \times 100$, to assess the error for the solution fields (with respect to FreeFEM results), the PDE residuals and boundary conditions. As can be seen in Tab.~\ref{tab:mse_fluid} and~\ref{tab:mse_solid}, the error is low, except for the residuals of the heat transfer equation for the fluid, which is localized in the vicinity of the origin, as already discussed.

\begin{table}[h]
\centering
\begin{tabular}{|l|c|}
    \hline
    Temperature               &  $6.6 \times 10^{-6}$ \\
    Velocity (u)              &  $1.0 \times 10^{-5}$ \\
    Velocity (v)              &  $7.3 \times 10^{-6}$ \\
    \hline
    No-slip             &  $3.4 \times 10^{-7}$ \\
    Inlet            &  $3.1 \times 10^{-8}$ \\
    Outlet           &  $2.9 \times 10^{-12}$\\
    Adiabaticity          &  $1.3 \times 10^{-5}$ \\
    \hline
    Heat eq.       &  $7.9 \times 10^{-3}$ \\
    Continuity eq. &  $7.1 \times 10^{-6}$ \\
    Momentum x eq. &  $1.2 \times 10^{-5}$ \\
    Momentum y eq. &  $7.0 \times 10^{-6}$ \\
    \hline
\end{tabular}
\caption{MSE for the physical fields, boundary conditions and PDE residuals for the fluid domain}
\label{tab:mse_fluid}
\end{table}

\begin{table}[h]
\centering
\begin{tabular}{|l|c|}
    \hline
    Temperature               &  $2.2 \times 10^{-6}$ \\
    \hline
    Dirichlet           &  $5.1 \times 10^{-8}$\\
    Adiabaticity          &  $1.4 \times 10^{-5}$ \\
    \hline
    Heat eq.       &  $1.1 \times 10^{-5}$ \\
    \hline
\end{tabular}
\caption{MSE for the physical fields, boundary conditions and PDE residuals for the solid domain}
\label{tab:mse_solid}
\end{table}

\clearpage

\section{Data assimilation}
\label{sec:assimilation}

\subsection{Method}

The ability of PINNs to handle both physical models and data, as well as their formulation as an optimization problem, make them especially well adapted to the resolution of data assimilation or inverse problems. In this context, some physical parameters of the problem $\gamma$ are unknown and one aims at estimating their values by minimizing the distance between the solution predicted $u(\gamma)$ and some observation data $u^\star$. Thus, for the system governed by Eq.~\ref{eq:pb_param}, for which $\gamma$ plays now the role of the unknown parameter, this can be expressed as an optimization problem contrained by PDEs:
\begin{equation}
    \min \mathcal{J}(\gamma) = \frac{1}{2} \| u - u^\star \|^2 \text{ s.t. } \frac{\partial{u}}{\partial{t}} + \mathcal{N}_\gamma(u,\frac{\partial{u}}{\partial{x}}) = 0
\end{equation}
Boundary and initial conditions are omitted here for the sake of readability. The classical way to solve such a problem is to embed the resolution of the physical system in an optimization loop, including an adjoint method to estimate the gradient of $\mathcal J$ with respect to $\gamma$ and perform a descent step. This approach is quite expensive because several PDE systems (state and adjoint) have to be solved in an iterative fashion.

Again, the specific formulation of PINNs allows to simplify the implementation and solve simultaneously the physics and the data assimilation problem, avoiding thus to introduce a second optimization loop. More precisely, the functional $\mathcal{J}(\gamma)$ is added in the loss function and its gradient is computed by automatic differentiation, as other loss terms. Then, the value of $\gamma$ is updated by a descent step, as the network parameters $\theta$. Therefore, the two optimization problems (determination of the network parameters and unknown physical parameters) are solved simulataneously.

In the following section, we apply this paradigm to solve a data assimilation problem involving turbulent flows.

\subsection{Application to backward facing step}

We consider the Reynolds-Averaged Navier-Stokes (RANS) equations for incompressible flows, with the Boussinesq hypothesis, denoted as $\mathcal{N}_{RANS}$:
\begin{equation}
    \frac{\partial{U}}{\partial x} + \frac{\partial{V}}{\partial y} = 0
\end{equation}
\begin{equation}
% \begin{split}
    \frac{\partial{U U}}{\partial x} + \frac{\partial{U V}}{\partial y}
    +\frac{1}{\rho} \frac{\partial{\tilde{P}}}{\partial x}  - \frac{\partial}{\partial x} \left ( (\nu + \nu_t) \frac{\partial U}{\partial x} \right )
     - \frac{\partial}{\partial y} \left ( (\nu + \nu_t) \frac{\partial U}{\partial y} \right ) = 0
% \end{split}
\end{equation}
\begin{equation}
% \begin{split}
    \frac{\partial{U V}}{\partial x} + \frac{\partial{V V}}{\partial y}
    +\frac{1}{\rho} \frac{\partial{\tilde{P}}}{\partial y}  - \frac{\partial}{\partial x} \left ( (\nu + \nu_t) \frac{\partial V}{\partial x} \right ) - \frac{\partial}{\partial y} \left ( (\nu + \nu_t) \frac{\partial V}{\partial y} \right ) = 0
% \end{split}
\end{equation}
with $\tilde{P} = P + \frac{2}{3} \rho K$, where $(U,V)$ represent the components of the mean velocity, $P$ the mean pressure, $\nu_t$ the turbulent kinematic viscosity and $K$ the turbulent kinetic energy. The resolution of $(U,V,P)$ necessitates to introduce a \emph{turbulence closure} to estimate $\nu_t$ and $K$. However, it is well known that turbulence models are limited in their range of application, therefore we aim at using PINNs to solve a data assimilation problem, based on observations of mean velocity profiles, to predict the turbulent viscosity field without introducing a turbulence model.

The flow over a backward facing step is considered as test-case, at moderate Reynolds number $R_e = 5100$, for which experimental and Direct Numerical Simulation (DNS) data are available. This configuration has already been studied using PINNs\cite{fluids8020043} and thus comparisons can be drawn regarding the results. The step height is $h=0.051$, the channel height is $6h$ and its length $23h$, the step being located at $x=3h$. We consider data collected at four sections located at $x=0$, $x=7$, $x=13$ and $x=22$, the mean velocity $U_i^\star$ at 22 observation points for each section being actually used.

\begin{figure*}[!ht]
\centering
\begin{tikzpicture}[x=1cm,y=0.5cm]

\node at (4,-4.5)  {descent step};
\draw[->,thick] (6,-5) -- (3,-5) -- (3,-3);

\draw[-,thick] (6,-5) -- (13,-5);
\node[mynode2] at (7,-5) {$\frac{\partial}{\partial{\theta}}$};

\draw[-,thick] (16.5,2) -- (17,0);
\draw[-,thick] (16.5,-1) -- (17,0);
\draw[-,thick] (16.5,-3) -- (17,0);
\draw[->,thick] (17,0) -- (17,-4);
\node at (16,-5)  {$\mathcal{L} = \mathcal{L}_{RANS} + \mathcal{L}_{W} + \mathcal{L}_{data}$} ;

\draw[-,thick] (5,1.5) -- (9,-3);
\draw[-,thick] (5,1.5) -- (9,-1);
\draw[-,thick] (5,0.5) -- (9,-1);

\node[right] at (9,-1)  {$\mathcal{L}_{W} = \frac{1}{N_{W}} \sum_{i=1}^{N_{W}} \left [ U(x_i^W,y_i^W) \right ]^2 + \left [ V(x_i^W,y_i^W) \right ]^2$};
\node[right] at (9,-3)  {$\mathcal{L}_{data} = \frac{1}{N_{data}} \sum_{i=1}^{N_{data}} \left [ U(x_i^{data},y_i^{data}) - U^\star_i \right ]^2$};

\draw[-,thick] (8,3) -- (9,2);
\draw[-,thick] (8,1) -- (9,2);
\node[right] at (9,2)  {$\mathcal{L}_{RANS} = \frac{1}{N_{RANS}} \sum_{i=1}^{N_{RANS}} \left [
\mathcal{N}_{RANS}(x_i,y_i)
 \right]^2$};

\draw[-,thick] (5,1.5) -- (7,3);
\draw[-,thick] (5,0.5) -- (7,3);
\draw[-,thick] (5,-0.5) -- (7,3);
\draw[-,thick] (5,-1.5) -- (7,3);
\draw[-,thick] (5,1.5) -- (7,1);
\draw[-,thick] (5,0.5) -- (7,1);
\draw[-,thick] (5,-0.5) -- (7,1);
\draw[-,thick] (5,-1.5) -- (7,1);

\draw[-,thick] (7,3) -- (8,3);
\draw[-,thick] (7,1) -- (8,1);

\node[mynode2] at (7,3)  {$\frac{\partial}{\partial x}$};
\node[mynode2] at (7,1)  {$\frac{\partial}{\partial y}$};
\node[mynode2] at (8,3)  {$\frac{\partial^2}{\partial x^2}$};
\node[mynode2] at (8,1)  {$\frac{\partial^2}{\partial y^2}$};

\node at (2.5,-2.5)  {$\theta$};
  \readlist\Nnod{2,5,5,5,4} % number of nodes per layer
  \foreachitem \N \in \Nnod{ % loop over layers
    \foreach \i [evaluate={\x=\Ncnt; \y=\N/2-\i+0.5; \prev=int(\Ncnt-1);}] in {1,...,\N}{ % loop over nodes
      \node[mynode] (N\Ncnt-\i) at (\x,\y) {};
      \ifnum\Ncnt>1 % connect to previous layer
        \foreach \j in {1,...,\Nnod[\prev]}{ % loop over nodes in previous layer
          \draw[-,thick] (N\prev-\j) -- (N\Ncnt-\i); % connect arrows directly
        }
      \fi % else: nothing to connect first layer
    }
  }
\node at (1,0.5)  {$x$};
\node at (1,-0.5)  {$y$};
\node at (5,1.5)  {$U$};
\node at (5,0.5)  {$V$};
\node at (5,-0.5)  {$\tilde{P}$};
\node at (5,-1.5)  {$\tilde{\nu}$};

\end{tikzpicture}
\caption{Data assimilation problem for the backward facing step case}
\label{fig:pinns_assimilation}
\end{figure*}

\begin{figure}[!ht]
\centering
\includegraphics[width=0.5\linewidth]{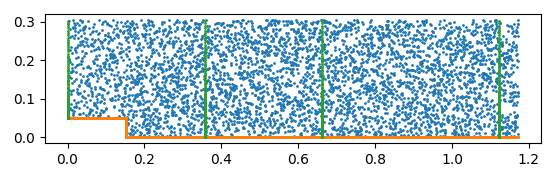}
\caption{Sampling and data sections}
\label{fig:bfs_sampling}
\end{figure}

The network here counts two inputs $(x,y)$, four outputs $(U,V,\tilde{P},\tilde{\nu})$ and five hidden layers with $\{ 8, 16, 32, 16, 8 \}$ neurons. Note that we solve the problem for $\tilde{P}$ directly because there is no way to distinguish $P$ from $K$. The network used is far smaller than the one employed in~\cite{fluids8020043}, which counts five layers of 128 neurons.
We enforce the positivity of the turbulent viscosity by setting $\nu_t = \tilde{\nu}^2$, which plays the role of the unknown field for the inverse problem. The training is achieved by minimizing the global loss function, composed of the MSE related to the RANS equations $\mathcal{L}_{RANS}$, the no-slip conditions at the wall $\mathcal{L}_{W}$ and the MSE related to the data $\mathcal{L}_{data}$, as depicted in Fig~\ref{fig:pinns_assimilation}.
We do not use additional terms to impose oulet boundary condition or velocity at the top boundary, contrary to~\cite{fluids8020043}, yielding a very simple implementation.
The loss terms related to RANS equations and no-slip conditions are based on samplings of size 8000 and 2750 points respectively, according to a Latin Hypercube distribution, as seen in Fig.~\ref{fig:bfs_sampling}.
The training is achieved using using ADAM algorithm (5,000 epochs) followed by BFGS algorithm (15,000 epochs).

\begin{figure}[!ht]
\centering
\includegraphics[width=0.5\linewidth]{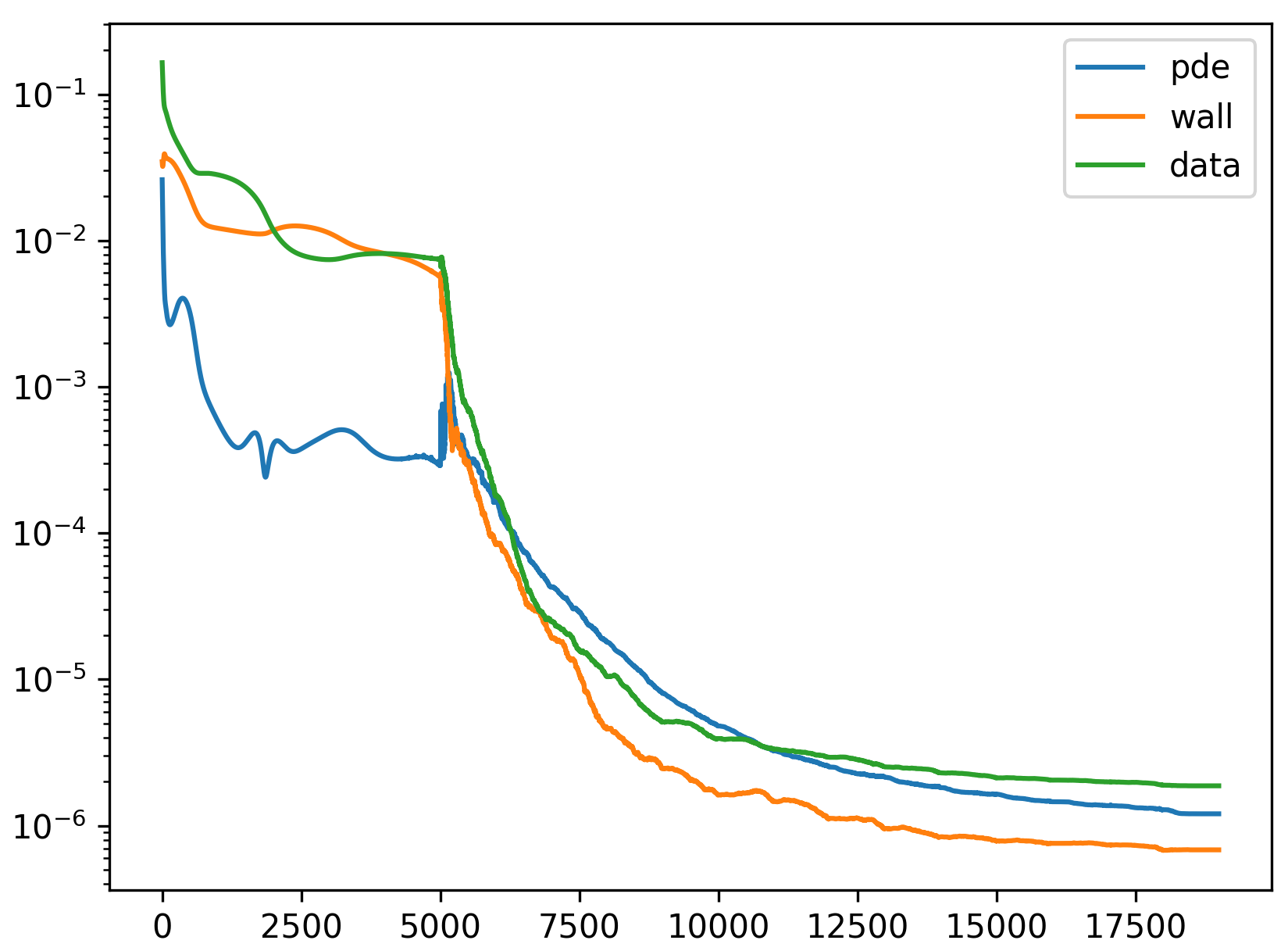}
\caption{Evolution of the loss terms}
\label{fig:bfs_losses}
\end{figure}

\begin{figure}[!ht]
\begin{minipage}{0.5\linewidth}
\centering
\includegraphics[width=0.9\linewidth]{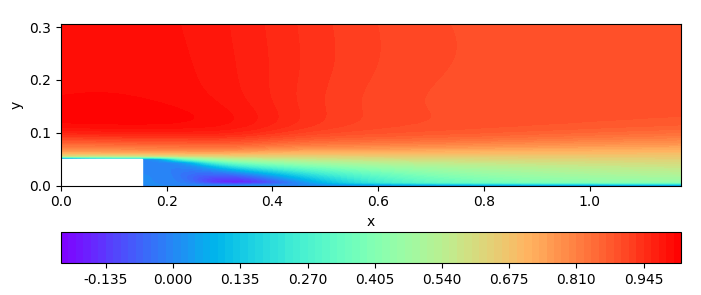}\\
\vspace{1mm}
\includegraphics[width=0.9\linewidth]{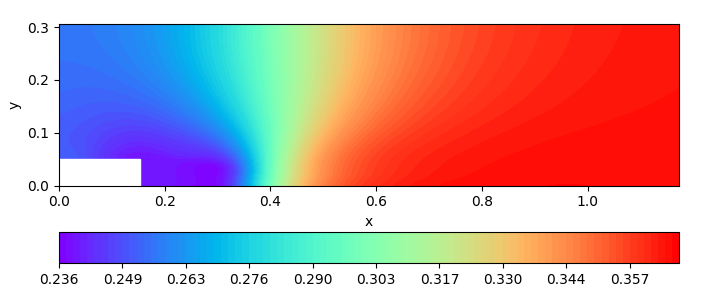}\\
\end{minipage}\hfill
\begin{minipage}{0.5\linewidth}
\centering
\includegraphics[width=0.9\linewidth]{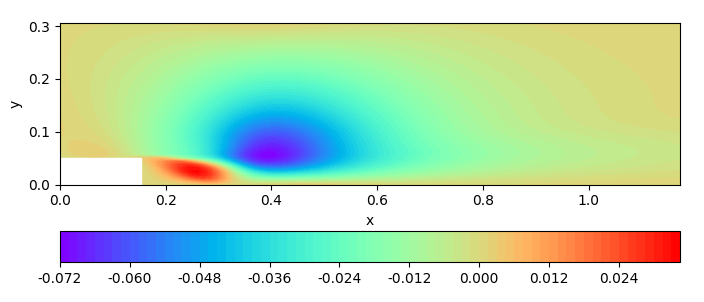}\\
\vspace{1mm}
\includegraphics[width=0.9\linewidth]{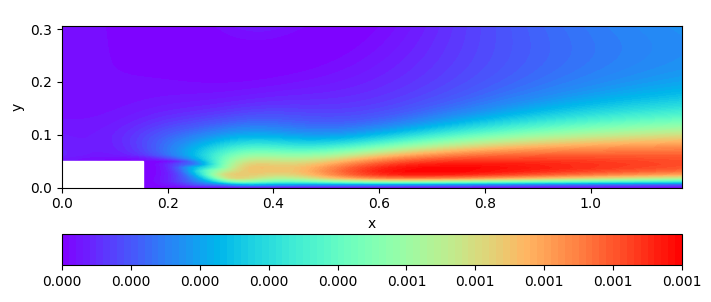}\\
\end{minipage}
\caption{Solution $U$ (top left), $V$ (top right), $\tilde{P}$ (bottom left) and $\nu_t$ (bottom right)}
\label{fig:bfs_sol}
\end{figure}

\begin{figure}[!ht]
\begin{minipage}{0.5\linewidth}
\centering
\includegraphics[width=0.9\linewidth]{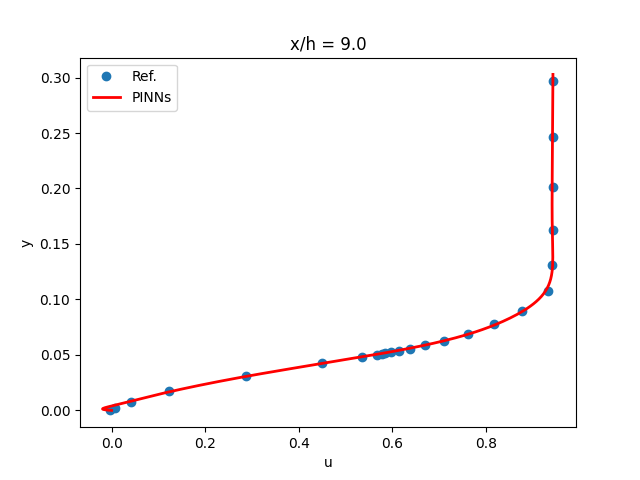}
\end{minipage}\hfill
\begin{minipage}{0.5\linewidth}
\centering
\includegraphics[width=0.9\linewidth]{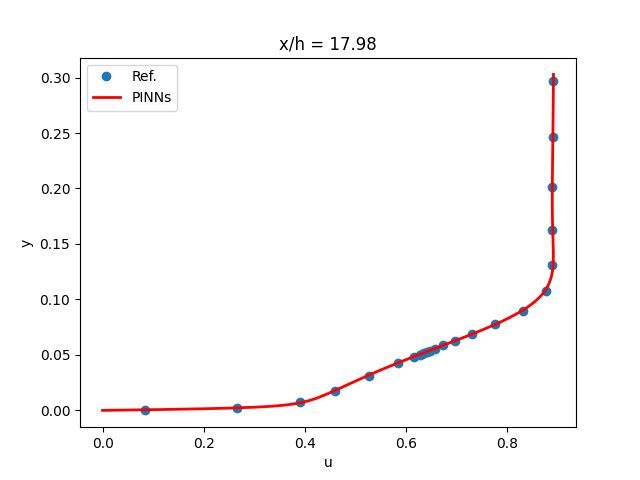}
\end{minipage}
\caption{Velocity profile at $x=9$ (left) and $x=18$ (right)}
\label{fig:bfs_vel}
\end{figure}

\begin{figure}[!ht]
\begin{minipage}{0.5\linewidth}
\centering
\includegraphics[width=0.9\linewidth]{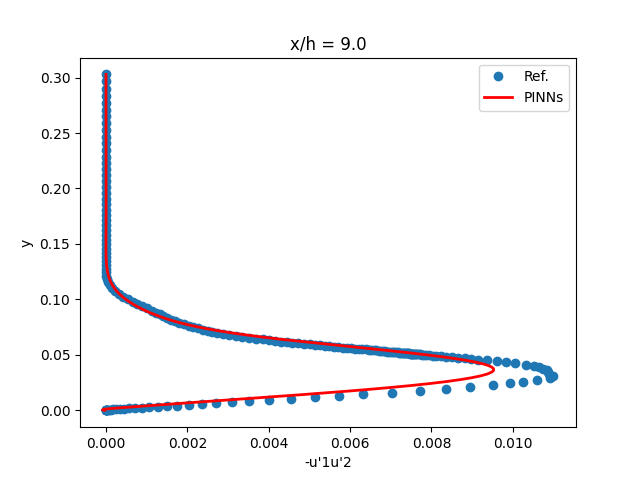}
\end{minipage}\hfill
\begin{minipage}{0.5\linewidth}
\centering
\includegraphics[width=0.9\linewidth]{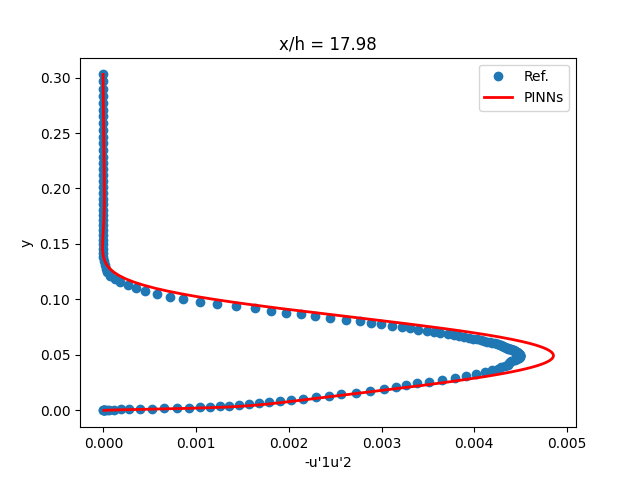}
\end{minipage}
\caption{Correlation profile at $x=9$ (left) and $x=18$ (right)}
\label{fig:bfs_corr}
\end{figure}

The convergence of the training process is shown in Fig.~\ref{fig:bfs_losses}. As seen, the use of a BFGS optimizer is critical for an efficient minimization of the loss function, ADAM algorithm being unable to achieve a significant reduction. The positivity enforcement of $\nu_t$ has been found necessary for a fast and reliable convergence. The solution fields obtained are depicted in Fig.~\ref{fig:bfs_sol}.
To assess the accuracy of these results, the DNS velocity profiles at sections $x=9$ and $x=18$, which are not included in the training set, are compared to those predicted by PINNs method in Figs.~\ref{fig:bfs_vel}.
As seen, these profiles are correctly inferred by the network. The correlation profiles $u'v'$, computed from the turbulent viscosity field thanks to the Boussinesq hypothesis are compared to the ones obtained from DNS (see Figs.~\ref{fig:bfs_corr}). A satisfactory agreement is observed, similar to those found in the literature using simulations based on standard tubulence closures.

\clearpage

\section{Discussion}

The above sections demonstrate the potentiality of PINNs to solve some non-classical PDE problems.
In this section, we aim at analysing more in depth the accuracy and robustness of the proposed method. In this perspective, the data assimilation problem presented in the previous section is considered as test-case and some statistical analyses are carried out. Some parameters of the method are modified and, for each case, ten trials are performed using different initializations of the network parameters $\theta$. Three error metrics are then estimated for the PDE residuals, the data fitting and the turbulent correlations inference, using sampling points and data not included in the training.

\begin{figure}[!ht]
\begin{minipage}{0.33\linewidth}
    \centering
    \hspace{3mm}\includegraphics[width=\textwidth]{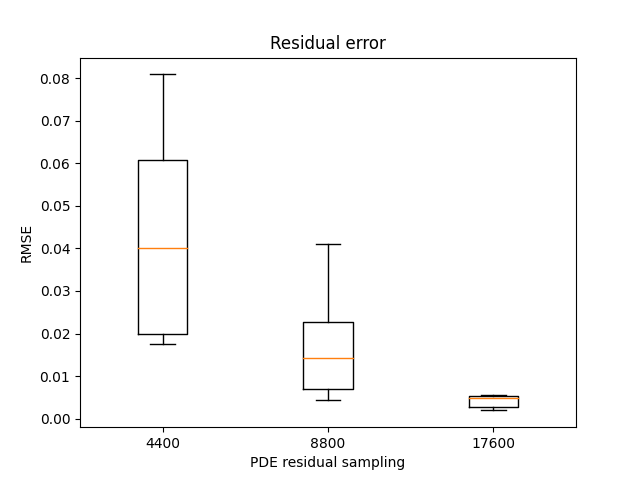}\\
\end{minipage}\hfill
\begin{minipage}{0.33\linewidth}
    \centering
    \hspace{3mm}\includegraphics[width=\textwidth]{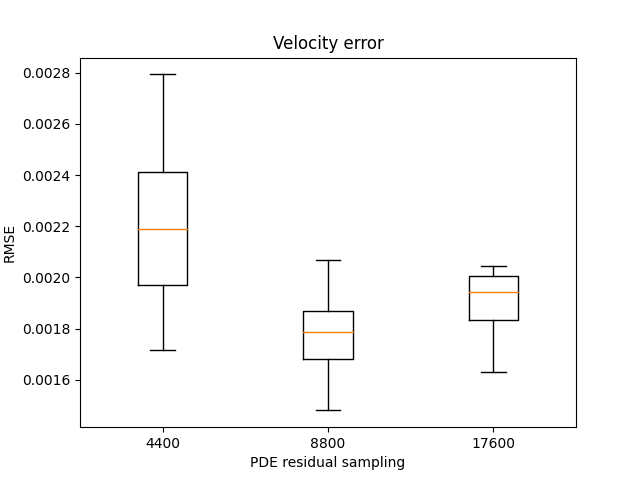}\\
\end{minipage}\hfill
\begin{minipage}{0.33\linewidth}
    \centering
    \hspace{3mm}\includegraphics[width=\textwidth]{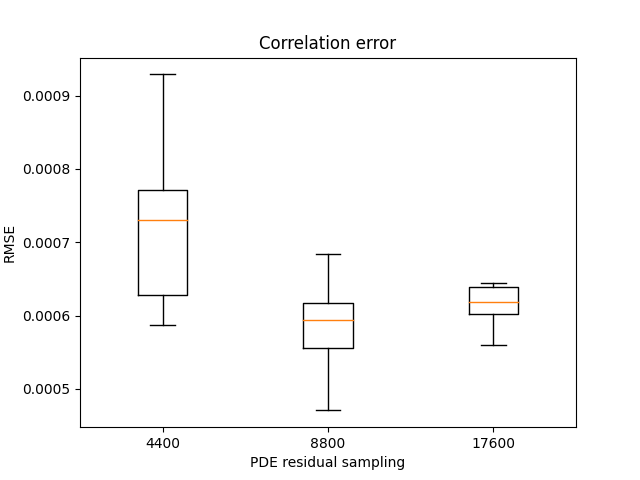}\\
\end{minipage}
\caption{Error w.r.t. sampling size}
\label{fig:box_pde_correlation}
\end{figure}

The results obtained using different sampling sizes are shown in Fig.~\ref{fig:box_pde_correlation}. As expected, an increase of the sampling size reduces the Root Mean-Square Error (RMSE) related to the residuals, in terms of average and variance. Note that the errors related to the data fitting and the  correlations inference are not so much impacted. Consequently, a very fine sampling is not necessary to solve the inverse problem.

\begin{figure}[!ht]
\begin{minipage}{0.33\linewidth}
    \centering
    \hspace{3mm}\includegraphics[width=\textwidth]{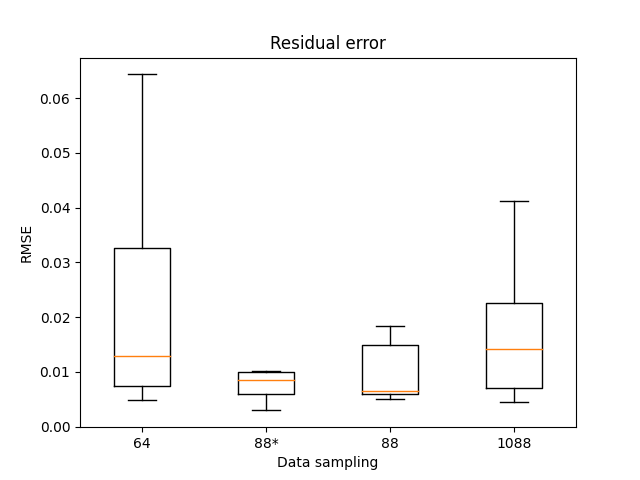}\\
\end{minipage}\hfill
\begin{minipage}{0.33\linewidth}
    \centering
    \hspace{3mm}\includegraphics[width=\textwidth]{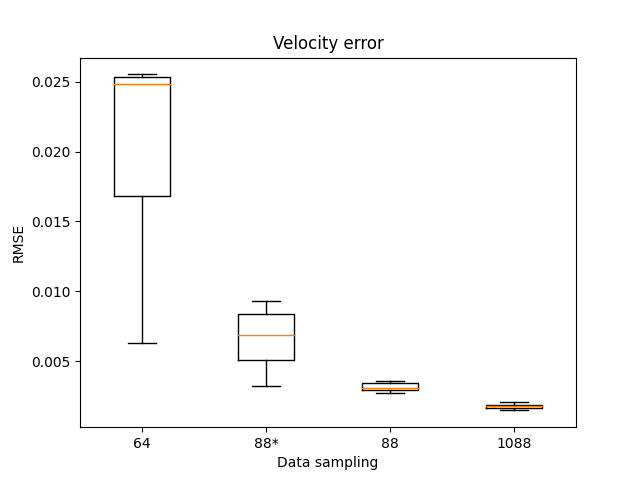}\\
\end{minipage}\hfill
\begin{minipage}{0.33\linewidth}
    \centering
    \hspace{3mm}\includegraphics[width=\textwidth]{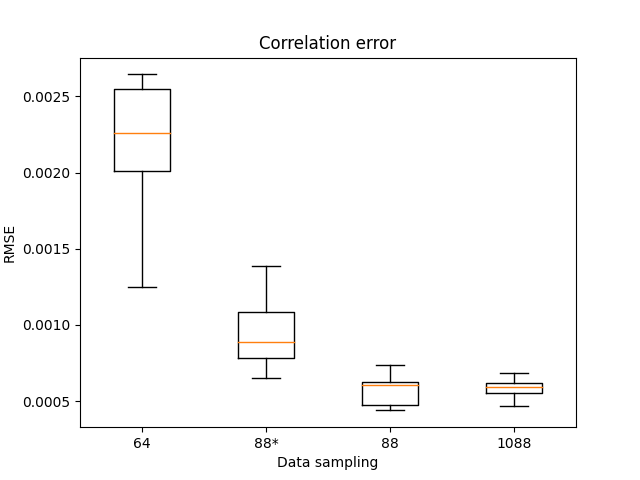}\\
\end{minipage}
\caption{Error w.r.t. sampling size}
\label{fig:box_data_correlation}
\end{figure}

The results obtained using different data sizes are shown in Fig.~\ref{fig:box_data_correlation}. As seen, increasing the data size yields a decrease of the data fitting error and correlations inference error. However, the improvement obtained using 1088 data points is rather small compared to the error obtained using only 88 points. The case with $88^\star$ points corresponds to a configuration where no data point is located in the recirculation area. As shown, a significant error increase is reported, which underlines the necessity to select carefully the data used for the training.

\begin{figure}[!ht]
\begin{minipage}{0.33\linewidth}
    \centering
    \hspace{3mm}\includegraphics[width=\textwidth]{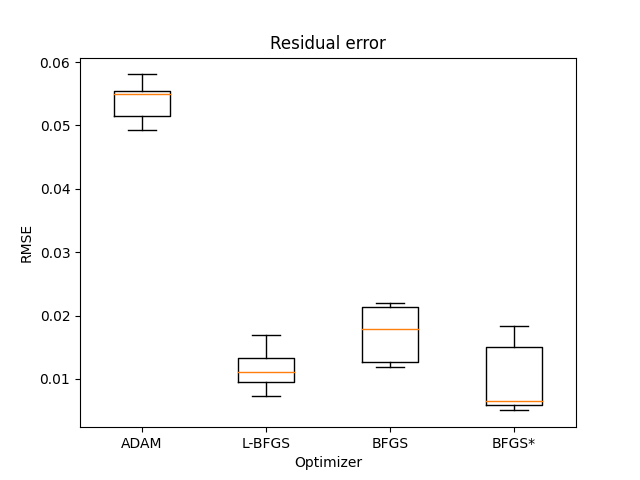}\\
\end{minipage}\hfill
\begin{minipage}{0.33\linewidth}
    \centering
    \hspace{3mm}\includegraphics[width=\textwidth]{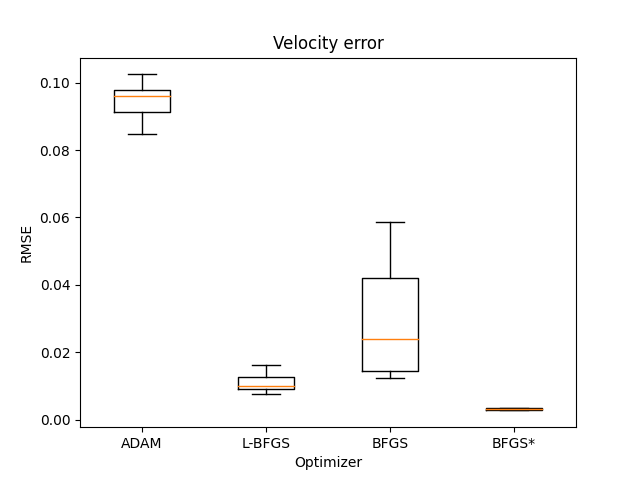}\\
\end{minipage}\hfill
\begin{minipage}{0.33\linewidth}
    \centering
    \hspace{3mm}\includegraphics[width=\textwidth]{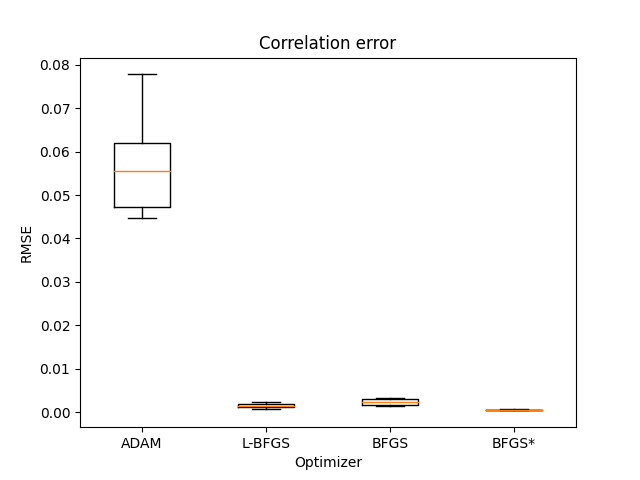}\\
\end{minipage}
\caption{Error w.r.t. optimization algorithm}
\label{fig:box_opt_correlation}
\end{figure}

The choice of the optimizer used to minimize the loss function is a critical parameter, as shown in Fig.~\ref{fig:box_opt_correlation}. Here, we compared the results obtained using ADAM only (20,000 epochs) against those using ADAM (5,000 epochs) followed by a quasi-Newton algorithm (15,000 epochs). BFGS and L-BFGS correspond to the implementations found in TensorFlow Probability 0.21 (Hager-Zhang line search), whereas BFGS$^\star$ corresponds to our own implementation based on Armijo-Goldstein line search. As seen, ADAM alone yields poor results and our implementation of BFGS performs far better than the one implemented in TensorFlow library. It shows that results can strongly depend on algorithmic details of the optimizer.

\begin{figure}[!ht]
\begin{minipage}{0.33\linewidth}
    \centering
    \hspace{3mm}\includegraphics[width=\textwidth]{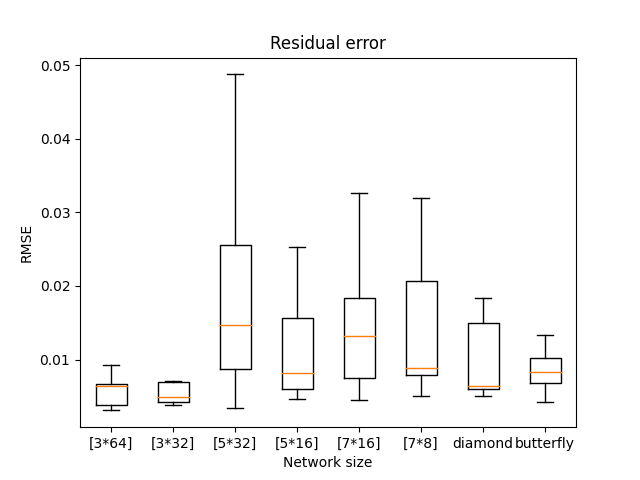}\\
\end{minipage}\hfill
\begin{minipage}{0.33\linewidth}
    \centering
    \hspace{3mm}\includegraphics[width=\textwidth]{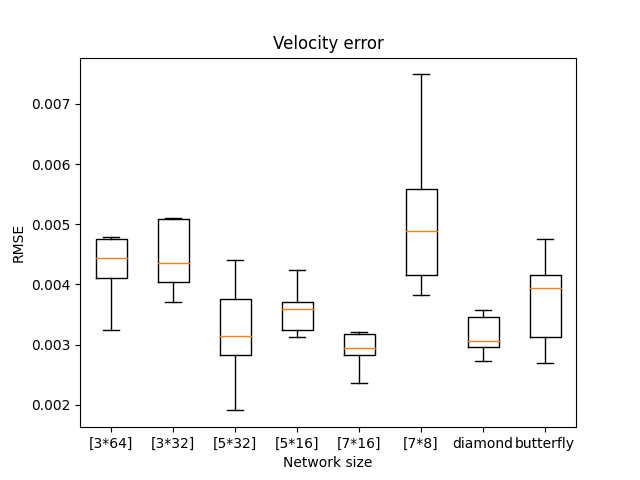}\\
\end{minipage}\hfill
\begin{minipage}{0.33\linewidth}
    \centering
    \hspace{3mm}\includegraphics[width=\textwidth]{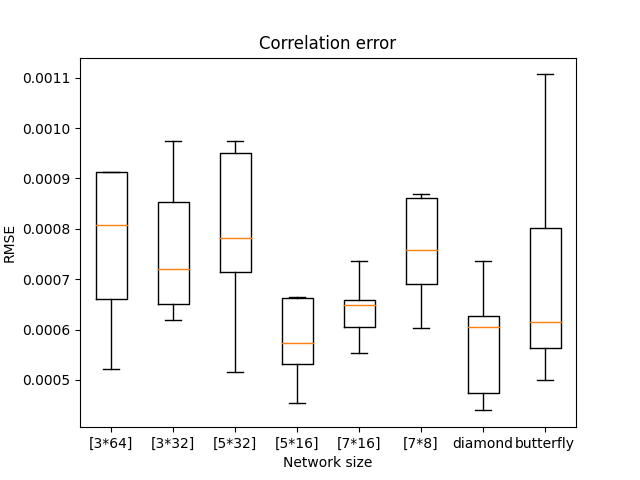}\\
\end{minipage}
\caption{Error w.r.t. network shape and size}
\label{fig:box_nn_correlation}
\end{figure}

The results obtained using different network shapes and sizes are shown in Fig.~\ref{fig:box_nn_correlation}. We test different numbers of layers, between 3 and 7, and different numbers of neurons per layer, between 8 and 64. The "diamond" corresponds to $\{8, 16, 32, 16, 8\}$ neurons and the "butterfly" to $\{32, 16, 8, 16, 32\}$ neurons (auto-encoder structure). It appears that satisfactory results can be obtained with all the configurations tested. A small number of layers provides better results regarding the residuals, but the resolution of the inverse problem performs better for deeper networks.

\begin{figure}[!ht]
\begin{minipage}{0.33\linewidth}
    \centering
    \hspace{3mm}\includegraphics[width=\textwidth]{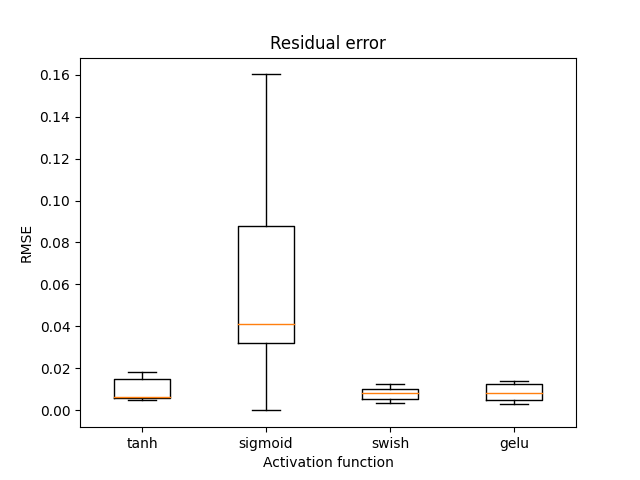}\\
\end{minipage}\hfill
\begin{minipage}{0.33\linewidth}
    \centering
    \hspace{3mm}\includegraphics[width=\textwidth]{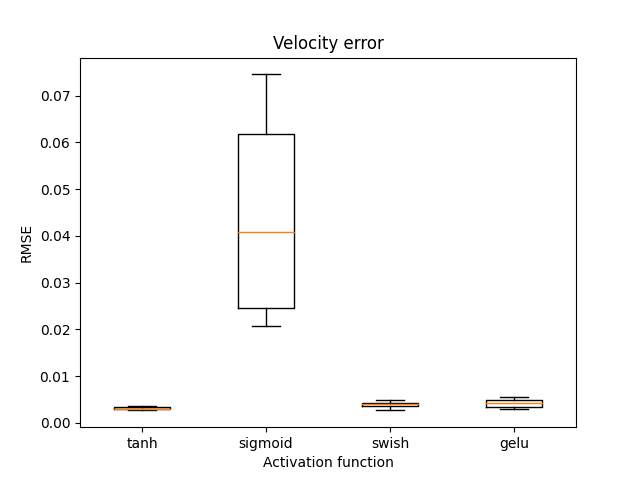}\\
\end{minipage}\hfill
\begin{minipage}{0.33\linewidth}
    \centering
    \hspace{3mm}\includegraphics[width=\textwidth]{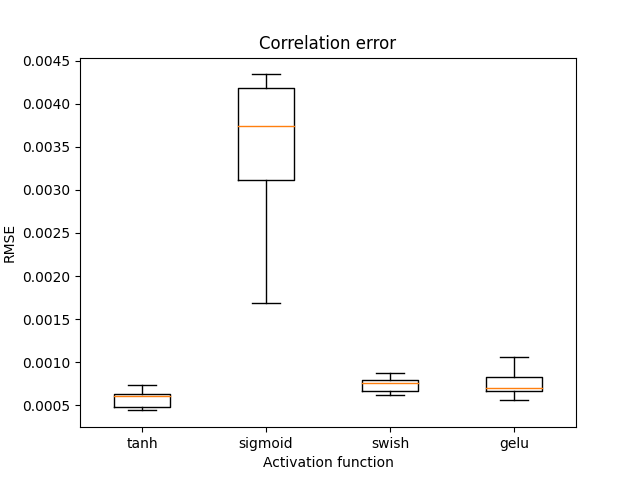}\\
\end{minipage}
\caption{Error w.r.t. activation functions}
\label{fig:box_activ_correlation}
\end{figure}

Finally, the impact of the choice of the activation function is shown in Fig.~\ref{fig:box_activ_correlation}. Four different functions are tested and satisfactory results are obtained with all except the sigmoid function. A possible reason could be related to the positivity of the function.

\clearpage

\section{Conclusion}

The three examples presented in this work show that the PINNs can be interesting as a complementary tool, for tasks for which conventional approaches are not very effective, such as  parametric modeling, multidisciplinary coupling or data assimilation.

However it should be underlined that, even if the implementation of PINNs is quite straightforward, several painful trials are usually necessary to tune the different numerical parameters and obtain finally satisfactory results. This puts in light the lack of understanding in the training, yielding a lack of control in terms of accuracy and robustness. This work indicates that the critical step is the optimization strategy employed for the loss minimization and, therefore, further studies have to be conducted regarding this topic.

% Bibliography
\bibliographystyle{plain}
\bibliography{references}

\end{document}